\documentclass[12pt]{article}
\begin{document}
\renewcommand{\thefootnote}{\fnsymbol{footnote}}
\newpage
\pagestyle{empty}
\setcounter{page}{0}
\newcommand{\dotimes}{\stackrel{}{\dot{\otimes}}}
\renewcommand{\thesection}{\arabic{section}}
\renewcommand{\theequation}{\thesection.\arabic{equation}}
\newcommand{\sect}[1]{\setcounter{equation}{0}\section{#1}}
\newfont{\twelvemsb}{msbm10 scaled\magstep1}
\newfont{\eightmsb}{msbm8}
\newfont{\sixmsb}{msbm6}
\newfam\msbfam
\textfont\msbfam=\twelvemsb
\scriptfont\msbfam=\eightmsb
\scriptscriptfont\msbfam=\sixmsb
\catcode`\@=11
\def\Bbb{\ifmmode\let\next\Bbb@\else
\def\next{\errmessage{Use \string\Bbb\space only in math mode}}\fi\next}
\def\Bbb@#1{{\Bbb@@{#1}}}
\def\Bbb@@#1{\fam\msbfam#1}
\newfont{\twelvegoth}{eufm10 scaled\magstep1}
\newfont{\tengoth}{eufm10}
\newfont{\eightgoth}{eufm8}
\newfont{\sixgoth}{eufm6}
\newfam\gothfam
\textfont\gothfam=\twelvegoth
\scriptfont\gothfam=\eightgoth
\scriptscriptfont\gothfam=\sixgoth
\def\frak{\frak@}
\def\frak@#1{{\fam\gothfam{{#1}}}}
\def\frak@@#1{\fam\gothfam#1}
\catcode`@=12
%
%
%
\def\CC{{\Bbb C}}
\def\NN{{\Bbb N}}
\def\QQ{{\Bbb Q}}
\def\RR{{\Bbb R}}
\def\ZZ{{\Bbb Z}}
\def\cA{{\cal A}}          \def\cB{{\cal B}}          \def\cC{{\cal C}}
\def\cD{{\cal D}}          \def\cE{{\cal E}}          \def\cF{{\cal F}}
\def\cG{{\cal G}}          \def\cH{{\cal H}}          \def\cI{{\cal I}}
\def\cJ{{\cal J}}          \def\cK{{\cal K}}          \def\cL{{\cal L}} 
\def\cM{{\cal M}}          \def\cN{{\cal N}}          \def\cO{{\cal O}}
\def\cP{{\cal P}}          \def\cQ{{\cal Q}}          \def\cR{{\cal R}} 
\def\cS{{\cal S}}          \def\cT{{\cal T}}          \def\cU{{\cal U}}
\def\cV{{\cal V}}          \def\cW{{\cal W}}          \def\cX{{\cal X}}
\def\cY{{\cal Y}}          \def\cZ{{\cal Z}}
\def\qed{\hfill \rule{5pt}{5pt}}
\def\id{\mbox{id}}
\def\ggo{{\frak g}_{\bar 0}}
\def\uqggo{\cU_q({\frak g}_{\bar 0})}
\def\uqggp{\cU_q({\frak g}_+)}
\def\half{\frac{1}{2}}
\def\btf{\bigtriangleup}
\newtheorem{lemma}{Lemma}
\newtheorem{prop}{Proposition}
\newtheorem{theo}{Theorem}
\newtheorem{Defi}{Definition}

\vfill
\vfill
\begin{center}

{\LARGE {\bf {\sf 
Modified braid equations, Baxterizations and noncommutative 
spaces for the quantum groups $GL_{q}(N), SO_{q}(N)$ and 
$Sp_{q}(N)$  
}}} \\[0.8cm]

{\large  
A. Chakrabarti$^{1}$\footnote{E-mail:
{\tt chakra@cpht.polytechnique.fr}} and R. Chakrabarti$^{2}$
\footnote{Permanent address:\,Department of Theoretical Physics, 
University of Madras, Guindy Campus, Chennai, 600 025, 
India,\,E-mail: {\tt ranabir@imsc.ernet.in}}
}

\begin{center}
{\em 
$^{1}$Centre de Physique Th\'eorique, Ecole Polytechnique,
91128 Palaiseau Cedex, France.\\
$^{2}$Institute of Mathematical Sciences, C. I. T. Campus,
Chennai, 600 113, India. \\[2.8cm]
}
\end{center}

\end{center}

\smallskip

\smallskip 

\smallskip

\smallskip

\smallskip

\smallskip 

\begin{abstract}
\noindent  Modified braid equations satisfied by generalized 
${\hat R}$ matrices ( for a {\em given} set of group relations 
obeyed by the elements of ${\sf T}$ matrices ) are constructed for 
$q$-deformed quantum groups $GL_q (N), SO_q (N)$ and $Sp_q (N)$  
with arbitrary values of $N$. The Baxterization of ${\hat R}$ 
matrices, treated as an aspect complementary to the 
{\em modification} of the braid equation, is obtained for all 
these cases in particularly elegant forms. A new class 
of braid matrices is discovered for the quantum groups $SO_{q}(N)$ 
and $Sp_{q}(N)$. The ${\hat R}$ matrices of this class, while being 
distinct from restrictions of the universal ${\hat{\cal R}}$ matrix 
to the corresponding vector representations, satisfy the standard 
braid equation. The modified braid equation and the Baxterization 
are obtained for this new class of ${\hat R}$ matrices. 
Diagonalization of the generalized ${\hat R}$ matrices is 
studied. The diagonalizers are obtained explicitly for some lower 
dimensional cases in a convenient way, giving directly the 
eigenvalues of the corresponding ${\hat R}$ matrices. Applications 
of such diagonalization are then studied in the context of 
associated covariantly quantized noncommutative spaces. 
\end{abstract}

\vfill
\newpage 

\pagestyle{plain}

\sect{Introduction}
Previously  one of us has studied `{\em modified} braid equation' 
(MBE) in the context of the quantum groups $GL_{p,q}(2), GL_{g,h}(2)$ 
and $GL_{q,h}(1|1)$ ( biparametric unitary, non-standard Jordanian 
and hybrid deformations respectively ) in \cite {C01}; and also for 
the orthogonal quantum group $SO_{q}(3)$ in \cite{C02}. The 
terminology is adopted from that of Gerstenhaber, Giaquinto and 
Schack \cite{GGS93,GG97}, who have studied a generalized class 
of deformations leading to MBE, where unlike the standard 
Yang-Baxter equation there are also inhomogeneous terms linear in 
tensored ${\hat R}$ matrices. These authors indicate the 
significance of and the interest in this equation. In the works 
\cite{C01,C02} the most general solutions of the quantum inverse 
scattering equation 
\begin{equation}
\hat{R} {\sf T_{1}} {\sf T_{2}} - {\sf T_{1}} {\sf T_{2}} 
\hat{R} = 0,\qquad {\sf T_{1}} = {\sf T} \otimes {\sf I},
\qquad {\sf T_{2}} = {\sf I} \otimes {\sf T}  
\label{eq:qis}
\end{equation}
for the relevant quantum groups were considered. Starting from a 
{\it given} set of group relations of the elements of the monodromy 
matrix ${\sf T}$, the most general ${\hat R}$ matrix satisfying
(\ref{eq:qis}) was constructed. It was observed that the standard
braid equation was modified for this generalized ${\hat R}$ matrix.
In this procedure the {\em conservation of the group laws postulated 
for the elements of the ${\sf T}$ matrix was maintained}. 

Here we construct, systematically and explicitly, the MBE for the 
quantized groups  $GL_{q}(N), SO_{q}(N)$ and $Sp_{q}(N)$ 
respectively. As explained in Sec. 1 of \cite{C02}, we 
will systematically exploit generalized spectral decomposition of 
the relevant ${\hat R}$ matrices in the vector representations. For 
an ${\hat R}$ matrix obeying the characteristic equation
\begin{equation}
\left({\hat R} - k_{1} {\sf I}\right)\, \left({\hat R} - k_{2} 
{\sf I}\right)\,\cdots \left({\hat R} - k_{p} {\sf I}\right) = 0,
\quad \left( k_{i} \neq k_{j}\,\,\,\hbox{if}\,\,i \neq j\,|\,( i, j)
\,=\,(1, 2,\cdots,p) \right),  
\label{eq:charact}
\end{equation}
the projectors in the eigenspaces of this ${\hat R}$ matrix read
\begin{equation}
P_{i} = \prod_{j \ne i} \frac{( {\hat R} -k_{j} {\sf I} )}
{( k_{i} - k_{j} )}, 
\label{eq:proj}
\end{equation} 
and satisfy the usual property
\begin{equation} 
P_{i}\,P_{j}\,=\,P_{i}\,\delta_{i\,j},\quad \sum_{i} P_{i}\,=\,{\sf I}.    
\label{eq:projrel}
\end{equation}
This orthonormalized set then provides the spectral decomposition.

Following \cite{RTF90} and \cite{I95} we first review the situation 
for the standard braid equation. In our discussion concerning 
braiding matrices, projectors and so on, we draw on \cite{RTF90}. In 
our study  of the Baxterization of the ${\hat R}$ matrices the 
analysis in \cite{I95} is particularly relevant. A large number of 
sources are cited in \cite{I95}, and we also refer to that list. The 
standard, {\em i.e.}`non-modified' in our context, braid equation 
reads
\begin{equation}
{\hat R}_{12}\,{\hat R}_{23}\,{\hat R}_{12}\, - \,
{\hat R}_{23}\,{\hat R}_{12}\,{\hat R}_{23}\,=\,0.
\label{eq:be}
\end{equation}
The spectral decomposition of the braid matrix of the quantum 
group $GL_{q}(N)$ with the conventional normalization \cite{I95}
is given by 
\begin{equation}
{\sf {\hat R}}\,=\,q P_{(+)} - q^{-1} P_{(-)},  
\label{eq:Aspcde}
\end{equation}
whereas braid matrices of the quantum groups $SO_{q}(N)$ 
and $Sp_{q}(N)\,( \hbox {where}\,N = 2n )$ in the usual 
normalization \cite{I95} may be written in a unified manner
\begin{equation}
{\sf{\hat R}} = q P_{(+)} - q^{-1} P_{(-)} + \varepsilon 
q^{\varepsilon - N} P_{(0)}.
\label{eq:BCDspcde}
\end{equation} 
In (\ref{eq:BCDspcde}) we have $\varepsilon = 1 (-1)$ for the 
quantum group $SO_{q}(N)\, \left(Sp_{q}(N)\right)$. The explicit 
expressions of the ${\hat R}$ matrices and the corresponding 
projectors for the above quantum groups are given in 
\cite{RTF90,I95}.

Now we proceed as follows. Maintaining the {\em same} projectors as 
in the standard braid equation, we generalize the ${\hat R}$ matrix
of the quantum group $GL_{q}(N)$ as
\begin{equation}
{\hat R} ( v ) = {\sf I}\,+\, v P_{(-)} = P_{(+)}\,+\, 
( 1 + v ) P_{(-)},
\label{eq:ARv}  
\end{equation}   
whereas the generalized ${\hat R}$ matrices for the quantum groups
$SO_{q}(N)$ and $Sp_{q}(N)$ read
\begin{equation}
{\hat R}( v, w ) = {\sf I} + v P_{(-)} + w P_{(0)}
= P_{(+)} + ( 1 + v ) P_{(-)} + ( 1 + w ) P_{(0)}.
\label{eq:BCDRvw}
\end{equation} 
Setting the `braid values' of the variables
\begin{equation}
v = - ( 1 + q^{\mp 2} ),\qquad w = - ( 1 - \varepsilon 
q^{\mp ( N + 1 - \varepsilon )} )
\label{eq:brdval}
\end{equation}
in (\ref{eq:ARv}) and (\ref{eq:BCDRvw}) we may recover the 
${\hat R}^{\pm 1}$ matrices satisfying the braid equation 
(\ref{eq:be}). Here and henceforth we will adopt the convention 
that for the braid matrices satisfying (\ref{eq:be}) we will not 
explicitly exhibit the corresponding values of the relevant spectral 
variables. We will throughout implement the normalizations used in 
(\ref{eq:ARv}) and (\ref{eq:BCDRvw}). For the purpose of later use 
we enlist here the braid matrices and their inverses according to
our normalization scheme. The ${\hat R}^{\pm 1}$ matrices of the 
quantum group $GL_q (N)$ read
\begin{equation}
{\hat R}^{\pm 1} = {\sf I} - ( 1 + q^{\mp 2} )\,P_{(-)},
\label{eq:RGLN}
\end{equation}
whereas the ${\hat R}^{\pm 1}$ matrices of the quantum groups
$SO_q (N)$ and $Sp_q (N)$ are given by
\begin{equation}
{\hat R}^{\pm 1} = {\sf I} - ( 1 + q^{\mp 2} )\,P_{(-)} 
- ( 1 - \varepsilon q^{\mp ( N + 1 - \varepsilon )} )\,P_{(0)}.
\label{eq:RRinvbrd}
\end{equation}
We observe that the braid matrices defined in (\ref{eq:RGLN})
and (\ref{eq:RRinvbrd}) following our normalization scheme
and the corresponding matrices  given respectively by 
(\ref{eq:Aspcde}) and (\ref{eq:BCDspcde}) as per the usual
normalization prescriptions, differ by an overall multiplicative 
factor:
\begin{equation}
{\sf{\hat R}} = q\,{\hat R}.
\label{eq:RqR}
\end{equation}     
The ${\hat R} ( v )$ and  ${\hat R} ( v, w )$, 
defined respectively in (\ref{eq:ARv}) and (\ref{eq:BCDRvw}), satisfy 
the following characteristic equations:
\begin{equation}
\left( {\hat R} ( v ) - {\sf I} \right)  
\left( {\hat R} ( v ) - ( 1 + v )\,{\sf I} \right) = 0,
\label{eq:Afact}
\end{equation}
\begin{equation}
\left( {\hat R} ( v, w ) - {\sf I} \right)
\left( {\hat R} ( v, w ) - ( 1 + v )\,{\sf I} \right) 
\left( {\hat R} ( v, w ) - ( 1 + w )\,{\sf I} \right) = 0.
\label{eq:BCDfact}
\end{equation} 
Using (\ref{eq:proj}) the relevant projectors can also be expressed 
directly as linear and quadratic functions of ${\hat R} ( v )$ and 
${\hat R } ( v, w )$ respectively. We also note that any ${\hat R}$
matrix satisfying ({\ref{eq:qis}) also necessarily satisfies     
\begin{equation}
f\left({\hat R}\right)\, {\sf T_{1}} {\sf T}_{2} - 
{\sf T_{1}} {\sf T_{2}} \,\,f\left({\hat R}\right) = 0,
\label{eq:fqis}
\end{equation}
where $f(x)$ is any well-behaved function. Due to the relations
(\ref{eq:Afact}) and (\ref{eq:BCDfact}) the operator $f({\hat R})$ 
reduces to a linear and a quadratic expression in the matrix 
${\hat R}$ in the respective cases. Thus for arbitrary values of 
the variables $v$ and $( v, w )$, our constructions in 
(\ref{eq:ARv}) and (\ref{eq:BCDRvw}) provide the {\em most general 
solutions} in the relevant examples.   
 
{\em Now comes the crucial question.} What modifications in the 
braid equation (\ref{eq:be}) are enacted as the variables $v$ and 
$( v, w )$ move away from the `braid values' given in 
(\ref{eq:brdval})? We are thus lead to our MBE for each case 
considered. As shown later in Sec. 2, the MBE for the quantum 
group $GL_{q}(N)$ reads 
\begin{equation}
{\hat R}_{12} ( v ) {\hat R}_{23} ( v ) {\hat R}_{12} ( v ) -
{\hat R}_{23} ( v ) {\hat R}_{12} ( v ) {\hat R}_{23} ( v ) =
c \left( {\hat R}_{12} ( v ) - {\hat R}_{23} ( v ) \right),
\label{eq:Ambe}
\end{equation}
whereas in the examples of the quantum groups $SO_{q}(N)$ and
$Sp_{q}(N)$  these equations have the form
\begin{eqnarray}
&&{\hat R}_{12} ( v,w ) {\hat R}_{23} ( v,w ) {\hat R}_{12} ( v,w ) 
- {\hat R}_{23} ( v,w ) {\hat R}_{12} ( v,w ) {\hat R}_{23} ( v,w ) 
\nonumber\\
&&\quad= c_{1} \left( {\hat R}_{12} ( v,w )
- {\hat R}_{23} ( v,w ) \right) 
+ c_{2} \left( {\hat R}_{12}^{-1} ( v,w ) 
- {\hat R}_{23}^{-1} ( v,w ) \right) \nonumber\\ 
&&\quad\phantom{=} + c_{3}  \left( {\hat R}_{12} ( v,w ) 
{\hat R}_{23}^{-1} ( v,w ) - {\hat R}_{23} ( v,w ) 
{\hat R}_{12}^{-1} ( v,w ) \right)\nonumber\\ 
&&\quad\phantom{=} - c_{3} \left( {\hat R}_{12}^{-1} ( v,w ) 
{\hat R}_{23} ( v,w ) - {\hat R}_{23}^{-1} ( v,w ) 
{\hat R}_{12} ( v,w ) \right). 
\label{eq:BCDmbe}
\end{eqnarray}
The coefficients $c$ and $( c_{1}, c_{2}, c_{3} )$ are given 
explicitly in Secs. 2 and 3 respectively. Inhomogeneous terms 
linearly depending on ${\hat R}$ matrix elements on the right hand 
side of the modified braid equation, as in (\ref{eq:Ambe}), were 
first considered \cite{GGS93, GG97} in the context of the quantum 
group $GL_{q}(N)$ with two projectors. There is, however, a sudden 
leap in complexity as the number of projectors increases by one as 
in the cases of quantum groups $SO_{q}(N)$ and $Sp_{q}(N)$. By 
restricting the parameters $v$ and $w$, we also select special 
cases in Secs. 3 and 4, where, for instance, $c_{3} = 0$ or 
$\{ c_{2} = 0, c_{3} = 0 \}$ holds. 

There is another aspect of our analysis. In \cite{C01} and 
\cite{C02}\,\,( particularly in Sec. 4 of \cite{C02} ) it was 
pointed out that the {\em modification} of the braid equation and 
the Baxterization of the ${\hat R}$ matrix are two complimentary 
facets of the generalized spectral decompositions in (\ref{eq:ARv}) 
and (\ref{eq:BCDRvw}). It is possible to proceed in one of the 
following two alternate directions. 

\smallskip   
\begin{itemize}
\item

The variables $v$ and $w$ are held fixed in each factor of 
the left, as in (\ref{eq:Ambe}) and (\ref{eq:BCDmbe}), and the 
inhomogeneous terms on the right are computed. This yields the MBE.

\item

The inhomogeneous terms on the right hand side of the braid 
equation may be constrained to be zero. This fixes the variables 
$v$ and $w$ in the appropriate ${\hat R}$ matrices in a particular 
fashion to be shown below. This provides the Baxterization of the 
${\hat R}$ matrices.  
\end{itemize}

The final results for the first possibility were presented above. 
Similarly we enlist below the final results for the other 
possibility. In Sec. 3.3 we study the additive form of 
Baxterization of the ${\hat R}$ matrices:
\begin{equation}
{\hat R}_{12} ( \theta )\,{\hat R}_{23} ( \theta + 
\theta^{\prime} )\,{\hat R}_{12} ( \theta^{\prime} )\,- 
{\hat R}_{23} ( \theta^{\prime} )
\,{\hat R}_{12} ( \theta + \theta^{\prime} )\,
{\hat R}_{23} ( \theta )\,=\,0,
\label{eq:adbax}
\end{equation}
where $q = {\hbox{exp}} (h)$. For this form of Baxterization the
variable $v ( \theta )$ reads
\begin{equation}
v (\theta) = \frac{\hbox{sinh} (h-\theta)}{\hbox{sinh} (h+\theta)} 
- 1,
\label{eq:vthta}
\end{equation}
for all the quantum groups studied here, namely $GL_{q}(N)$,
$SO_{q}(N)$ and $Sp_{q}(N)$. The variable $w ( \theta )$, appearing
for the quantum groups $SO_{q}(N)$ and $Sp_{q}(N)$, assumes two 
alternate forms. For the orthogonal quantum group $SO_{q}(N)$
it is given by
\begin{equation}
w ( \theta ) = \frac{\hbox{cosh}\left( \frac{N}{2} h - \theta 
\right)}{\hbox{cosh} \left( \frac{N}{2} h + \theta \right)} - 1
\label{eq:wSO1}
\end{equation}
or
\begin{equation}
w ( \theta ) = \frac{\hbox{sinh}\left( 
\left( \frac{N}{2} - 1 \right) h - \theta \right)\,\hbox{sinh} 
( h - \theta )}{\hbox{sinh} \left(
\left( \frac{N}{2} - 1 \right) h + \theta \right)\,\hbox{sinh}
( h + \theta )}\,-\,1.
\label{eq:wSO2}
\end{equation}  
For the symplectic quantum group $Sp_{q}(N)$, where $N = 2n$, the 
variable $w ( \theta )$ assumes the form
\begin{equation}
w ( \theta ) = \frac{\hbox{sinh} \left( \left( n + 1 \right) h - 
\theta \right)}{\hbox{sinh} \left( \left( n + 1 \right) h + \theta 
\right)} - 1
\label{eq:wSp1}
\end{equation}
or
\begin{equation}
w ( \theta ) = \frac{\hbox{cosh} \left( n h - \theta \right)\,
\hbox{sinh} ( h - \theta )}{\hbox{cosh} \left( n h + \theta \right)
\hbox{sinh} ( h + \theta )} - 1.
\label{eq:wSp2}
\end{equation}
In each case our parametrization and normalization assure the 
validity of the constraint: 
\begin{equation}
{\hat R} ( \theta )\, {\hat R} ( - \theta ) = {\sf I}.
\label{eq:unit}
\end{equation}

It has been implicitly assumed above that $v \neq 0, w \neq 0$.
But the parametrization (\ref{eq:BCDRvw}) insistently points at
the following question. What happens for the special choices 
$( v \neq 0, w = 0 )$ and $( v = 0, w \neq 0 )$? Certain 
properties of the projectors $P_{(-)}$ and $P_{(0)}$ being 
different, there is no reciprocal symmetry concerning the two 
above choices. For $w = 0$, there is no nontrivial solution. But 
for $v = 0$ we discuss {\em a hitherto unnoticed new class of 
solutions} studied in detail in Sec. 4. Even in the limit $q = 1$,
this class of solutions {\em remains nontrivial}. For the purpose 
of comparison with the preceding results, here we just exhibit the 
Baxterized values of the spectral variables corresponding to 
these solutions for the quantum groups $SO_{q}(N)$ and
$Sp_q (N)$:
\begin{equation}
v = 0,\,\,\, w ( \theta ) = \frac{\hbox{sinh} ( \eta - \theta )}
{\hbox{sinh} ( \eta + \theta )} - 1,\,\,\,\hbox{where}\quad
\hbox{tanh}\,\eta = \left( 1 - 4 ( 1 + \varepsilon 
[ N - \varepsilon ] )^{-2} \right)^{1/2}.
\label{eq:v0wbax}
\end{equation}
In (\ref{eq:v0wbax}) the quantity under the radical sign is 
positive for $N > 2$, and hence its square root is real.
 
In Sec. 5 we construct matrices diagonalizing the generalized
braid operators ${\hat R} (v)$ and ${\hat R} ( v, w )$. Explicit 
results are presented for the quantum groups $GL_{q}( 2 ), 
SO_{q} ( 3 )$ and $SO_{q} ( 4 )$. The key result valid for the 
cases studied is that the indefinite parameters in the 
diagonalizers may be chosen to ensure {\em the mutual 
orthogonality of their rows}. This remarkably helpful property 
allows us to effortlessly obtain the eigenvectors of the 
generalized braiding matrices ${\hat R} ( v )$ and 
${\hat R} ( v, w )$. These can be of interest 
in related statistical models.

Also we show, in Sec. 6, how such diagonalizations may be 
exploited in the description  of associated noncommutative spaces. 
One principal objective in introducing our generalized spectral 
decompositions has been \cite{C01,C02} the exploration of the 
roles of the variables $( v, w )$ in the instances, where the 
generalized braid matrices ${\hat R} ( v )$ and ${\hat R} ( v, w )$ 
are used to construct the relevant noncommutative spaces. The 
constraints due to Leibnitz rule and the covariance 
properties - discussed, for instance, in \cite{WZ90,M99,CFM00} and 
a large number of sources cited therein - are presented in Sec. 6 
after incorporating our spectral variables.

But in this paper our study remains essentially limited to showing
how the the diagonalization of Sec. 5 can help in better 
understanding certain features of such spaces. We hope to develop 
other aspects elsewhere. Remarks on various features of the 
results obtained and other perspectives are presented in Sec. 7.      
    
\sect{Modified braid equation and Baxterization for the quantum
group $GL_{q}(N)$}
In the vector representation the braid matrix ${\hat R}$ of the 
quantum group $GL_{q}(N)$ has two orthogonal projectors. The 
spectral decomposition of its generalized ${\hat R} ( v )$ matrix 
has been defined in (\ref{eq:ARv}). Adopting the definitions 
\begin{equation}
X_{1} = P_{(-)} \otimes {\sf I},\,\,X_{2} = {\sf I} \otimes P_{(-)},
\,\,S_{1} = \left(\,X_{1} - X_{2}\,\right),\,\,
T_{1} = \left(\,X_{1} X_{2} X_{1} - X_{2} X_{1} X_{2}\,\right),
\label{eq:Adef}
\end{equation} 
we write
\begin{equation}
{\hat R}_{12}(v) {\hat R}_{23}(v^{\prime}) 
{\hat R}_{12}(v^{\prime\prime}) -  {\hat R}_{23}(v^{\prime\prime}) 
{\hat R}_{12}(v^{\prime}) {\hat R}_{23}(v) = ( v + v^{\prime\prime} 
+ v v^{\prime\prime} -v^{\prime} ) S_{1} + v v^{\prime} 
v^{\prime\prime} T_{1}.
\label{eq:Av123}
\end{equation}
Equating the variables to their braid values $v = v^{\prime} = 
v^{\prime\prime} = - ( 1 + q^{- 2} )$ the right hand side of 
(\ref{eq:Av123}) vanishes as in this limit the braid matrix 
${\hat R}$ in 
(\ref{eq:RGLN}) satisfying the braid equation (\ref{eq:be}) is 
obtained. This restricts 
\begin{equation}
T_{1} = \frac{S_{1}}{[2]^2},\qquad \hbox{where}\,\, 
[x] = \frac{q^{x} - q^{- x}}{q - q^{-1}}.
\label{eq:Vrel}
\end{equation}  
Hence the right hand side of (\ref{eq:Av123}) reads
\begin{equation}
\left( ( v + v^{\prime\prime} + v v^{\prime\prime} - v^{\prime} )\,
+\,\frac{v v^{\prime} v^{\prime\prime}}{[2]^2}\right)\,S_{1},
\label{eq:Avrhs}
\end{equation}  
where following (\ref{eq:ARv}) and the related discussions we obtain
\begin{equation}
S_{1} = v^{-1} \left( {\hat R}_{12} (v) - {\hat R}_{23} (v) \right) 
= - q\,[2]^{-1}\,\left( {\hat R}_{12} - {\hat R}_{23} \right).
\label{eq:S1val}
\end{equation}
Setting $v = v^{\prime} = v^{\prime\prime}$ we obtain the MBE 
obeyed by the generalized braid matrix ${\hat R} (v)$ of the quantum 
group $GL_{q}(N)$: 
\begin{eqnarray}   
&&{\hat R}_{12}(v) {\hat R}_{23}(v) {\hat R}_{12}(v) - 
{\hat R}_{23}(v) {\hat R}_{12}(v) {\hat R}_{23}(v)\nonumber\\
&&\qquad= \left( 1 + v + [2]^{-2}\, v^2 \right) 
\left( {\hat R}_{12}(v) - {\hat R}_{23}(v) \right).
\label{eq:AMBEv}  
\end{eqnarray}  
Comparing the above MBE with (\ref{eq:Ambe}) we obtain
\begin{equation}
c = 1 + v + [2]^{-2}\, v^2.
\label{eq:cvalue}
\end{equation} 

For the Baxterization of the ${\hat R} ( v )$ matrix we set
\begin{equation}
v = v(x), v^{\prime\prime} = v(y), v^{\prime} = v(xy)
\label{eq:vbax}
\end{equation}
and denote ${\hat R} (v) \equiv {\hat R} (x)$. Now for the right 
hand side of (\ref{eq:Av123}) to vanish, we must have
\begin{equation}  
v(xy) = \frac{v(x) + v(y) + v(x) v(y)}{1 - [2]^{-2} v(x) v(y)}.
\label{eq:vbaxrel}
\end{equation}
This functional relation is solved in a more general form in 
Sec. 3. Here we present the final result concerning the function 
$v ( x )$ as follows. The solution satisfying the constraint 
\begin{equation} 
{\hat R} (x) {\hat R} \left(x^{-1} \right) = 1
\label{eq:unitary}
\end{equation}
reads
\begin{equation}
v(x) = \frac{q x^{-1} - q^{-1} x}{q x - q^{-1} x^{-1}} - 1.
\label{eq:univx}
\end{equation}
Setting $x = \hbox{exp} (\theta)$ and $q = \hbox{exp} (h)$ in
(\ref{eq:univx}), we 
obtain the functional structure in (\ref{eq:vthta}). 
For the choice $y = \hbox{exp} ( \theta^{\prime} )$ and the above 
value of $v ( \theta )$ the braid equation takes the form
given in (\ref{eq:adbax})..
Combining (\ref{eq:ARv}), (\ref{eq:RGLN}) and (\ref{eq:univx}) 
we may express 
${\hat R} (x)$ as
\begin{equation}
\hat{R} (x) = \frac{q x \hat{R} - q^{-1} x^{-1} \hat{R}^{-1}}
{q x - q^{-1} x^{-1}}.
\label{eq:Rx}
\end{equation} 
Suitably changing the normalizations, namely, observing 
(\ref{eq:RqR}) and setting
\begin{equation}
{\sf {\hat R}} ( x ) = \left( q x^{- 1} - q^{-1} x \right) 
{\hat R} ( x^{-1} ), 
\label{Rredef}
\end{equation}
we obtain
\begin{equation}
{\sf {\hat R}} ( x ) = x^{-1}\,{\sf{\hat R}} - 
x\,{\sf{\hat R}}^{-1}.
\label{Rusual}
\end{equation}
In this form the Baxterization of the ${\hat R}$ matrix of the 
quantum group $GL_q (N)$ is often presented  \cite{I95}. We 
have preferred (\ref{eq:Rx}) to achieve uniform utilization of the 
same functional equation in Secs. 2. 3 and 4. This also ensures 
{\em one uniform normalization prescription}, whether the 
${\hat R}$ matrix is Baxterized or not, by fixing the top left 
element (row 1, column 1) to be unity. This condition is evidently 
satisfied by (\ref{eq:Rx}). 

\sect{Modified braid equation and Baxterization for the quantum 
groups $SO_{q}(N)$ and $Sp_{q}(N)$} 
{\large{\bf 3.1 Reduction of trilinear terms}}\\ 
\medskip
\noindent For arbitrary parameters $(v, w)$ the generalized 
${\hat R} ( v, w )$ matrix is given by (\ref{eq:BCDRvw}). In 
addition to quantities already introduced in (\ref{eq:Adef}), here
we further define the following objects:
\begin{eqnarray}
&&Y_{1} = P_{(0)} \otimes {\sf I},\qquad
Y_{2} = {\sf I} \otimes P_{(0)},\qquad S_{2} = ( Y_{1} -Y_{2} ),
\nonumber\\
&&J_{1} = ( X_{1} Y_{2} - Y_{1} X_{2} ),\qquad
J_{2} = ( Y_{2} X_{1} - X_{2} Y_{1} ),\nonumber\\
&&K_{1} = ( X_{1} X_{2} Y_{1} - Y_{2} X_{1} X_{2} ),\,\,
K_{2} = ( X_{1} Y_{2} X_{1} - X_{2} Y_{1} X_{2} ),\,\,
K_{3} = ( Y_{1} X_{2} X_{1} - X_{2} X_{1} Y_{2} ),\nonumber\\
&&L_{1} = ( Y_{1} Y_{2} X_{1} - X_{2} Y_{1} Y_{2} ),\,\,
L_{2} = ( Y_{1} X_{2} Y_{1} - Y_{2} X_{1} Y_{2} ),\,\,
L_{3} = ( X_{1} Y_{2} Y_{1} - Y_{2} Y_{1} X_{2} ),\nonumber\\
&&S_{3} = ( K_{1} + K_{2} + K_{3} ),\,\, 
S_{4} = ( L_{1} + L_{2} + L_{3} ),\,\,
T_{2} = ( Y_{1} Y_{2} Y_{1} - Y_{2} Y_{1} Y_{2} ).
\label{eq:SOdef}
\end{eqnarray}
In terms of these quantities we obtain
\begin{eqnarray}
&&{\hat R}_{12} (v, w) {\hat R}_{23} (v^{\prime}, w^{\prime}) 
{\hat R}_{12} (v^{\prime\prime}, w^{\prime\prime}) - 
{\hat R}_{23} (v^{\prime\prime}, w^{\prime\prime})  
{\hat R}_{12} (v^{\prime}, w^{\prime})  
{\hat R}_{23} (v, w)\nonumber\\  
&&\quad= ( v + v^{\prime\prime} + v v^{\prime\prime} - 
v^{\prime} ) S_{1} + ( w + w^{\prime\prime} + w w^{\prime\prime} - 
w^{\prime} ) S_{2}\nonumber\\
&&\phantom{\quad=}+ ( v w^{\prime} - v^{\prime} w ) J_{1} 
+ ( w^{\prime} v^{\prime\prime} - w^{\prime\prime} v^{\prime} ) 
J_{2}\nonumber\\
&&\phantom{\quad=}+ ( v v^{\prime} w^{\prime\prime} ) K_{1}
+ ( v w^{\prime} v^{\prime\prime} ) K_{2}
+ ( w v^{\prime} v^{\prime\prime} ) K_{3}\nonumber\\
&&\phantom{\quad=}+ ( w w^{\prime} v^{\prime\prime} ) L_{1}
+ ( w v^{\prime} w^{\prime\prime} ) L_{2}
+ ( v w^{\prime} w^{\prime\prime} ) L_{3}\nonumber\\
&&\phantom{\quad=}+ ( v v^{\prime} v^{\prime\prime} ) T_{1}
+ ( w w^{\prime} w^{\prime\prime} ) T_{2}.
\label{eq:SORvw}
\end{eqnarray}
We will discuss below how to express the trilinear combinations 
$( K_1, K_2, K_3, L_1, L_2, L_3, T_1, T_2 )$ in terms of the 
linear $( S_1, S_2 )$ and the bilinear $( J_1, J_2 )$ constructs.
For the `braid values' of the parameters
\begin{equation}
v = v^{\prime} = v^{\prime\prime} = - ( 1 + q^{-2} ),\qquad
w = w^{\prime} = w^{\prime\prime} = 
- \bigl( 1 - \varepsilon q^{- ( N + 1 - \varepsilon )} \bigr)
\label{eq:brdvw}
\end{equation}
the right hand side of (\ref{eq:SORvw}) vanishes. As noted before,
when the spectral variables assume their `braid values', the 
generalized braid matrix ${\hat R} ( v, w )$, while satisfying
(\ref{eq:be}), reduces to the ${\hat R}$ matrix given in 
(\ref{eq:RRinvbrd}). The characteristic equation 
(\ref{eq:charact}) for the ${\hat R}$ matrix now reads
\begin{equation}
\left( {\hat R} - {\sf I} \right) 
\left( {\hat R} + q^{-2}\,{\sf I} \right) 
\left( {\hat R} - \varepsilon q^{- ( N + 1 - \varepsilon )}
\,{\sf I} \right)\,=\,0.
\label{eq:SOSpchr}
\end{equation}
The projectors may be extracted {}from (\ref{eq:SOSpchr}) as 
quadratic expression in the matrix ${\hat R}$ {\em \`a\, la} 
(\ref{eq:proj}). For the choice
\begin{equation}
v = v^{\prime} = v^{\prime\prime} = - ( 1 + q^{2} ),\qquad
w = w^{\prime} = w^{\prime\prime} = 
- \bigl( 1 - \varepsilon q^{N + 1 - \varepsilon} \bigr)
\label{eq:invbrdvw}
\end{equation}
the inverse of the braiding matrix, namely ${\hat R}^{-1}$ is 
obtained. The braid matrix and its inverse ${\hat R}^{\pm 1}$
obtained following the choices (\ref{eq:brdvw}) and 
(\ref{eq:invbrdvw}) have been listed in (\ref{eq:RRinvbrd}). For 
these matrices, as mentioned before, we will suppress the 
particular values of the spectral variables. For future use here 
we note that the inverse of the generalized braid matrix 
(\ref{eq:BCDRvw}) is given by
\begin{equation}
\left({\hat R} ( v, w ) \right)^{-1} = {\sf I} - v ( 1 + v )^{-1} 
P_{(-)} - w ( 1 + w )^{-1} P_{(0)}.
\label{eq:invRvw}
\end{equation}  

The braid equation (\ref{eq:be}) may be utilized to yield the
well-known constraints
\begin{equation}
f \left( \hat{R}_{12}^{\pm 1} \right) \,\,\hat{R}_{23}^{\pm 1} 
\hat{R}_{12}^{\pm 1} = \hat{R}_{23}^{\pm 1} \hat{R}_{12}^{\pm 1} 
\,\,f \left( \hat{R}_{23}^{\pm 1} \right),\qquad
\hat{R}_{12}^{\pm 1} \hat{R}_{23}^{\pm 1}\,\, 
f \left( \hat{R}_{12}^{\pm 1} \right) = 
f \left( \hat{R}_{23}^{\pm 1} \right)\,\,\hat{R}_{12}^{\pm 1} 
\hat{R}_{23}^{\pm 1}, 
\label{eq:becons}
\end{equation}
where $f(x)$ is any well-behaved function of $x$. The projectors 
$P_{(-)}$ and $P_{(0)}$ are quadratic functions of the matrix 
$\hat{R}$. Choosing these particular functions in the constraints 
(\ref{eq:becons}), and utilizing the definitions (\ref{eq:Adef}) 
and (\ref{eq:SOdef}) we obtain
\begin{equation}
X_{1} A_{(\pm)} = A_{(\pm)} X_{2},\quad 
Y_{1} A_{(\pm)} = A_{(\pm)} Y_{2},\quad 
B_{(\pm)} X_{1} = X_{2} B_{(\pm)},\quad
B_{(\pm)} Y_{1} = Y_{2} B_{(\pm)},
\label{eq:XYcon}
\end{equation}
where the bilinear elements read
\begin{eqnarray}
&&A_{(\pm)} = \hat{R}_{23}^{\pm 1} \hat{R}_{12}^{\pm 1}\nonumber\\ 
&&\phantom{A_{(\pm)}} = {\sf I} - \bigl( 1 + q^{\mp 2}\bigr)
( X_1 + X_2 ) - 
\left( 1 - \varepsilon q^{\mp ( N + 1 - \varepsilon)} \right)
( Y_1 + Y_2 ) + \left( 1 + q^{\mp 2} \right)^2 X_2 X_1\nonumber\\ 
&&\phantom{A_{(\pm)} =} + ( 1 + q^{\mp 2} ) 
\left( 1 - \varepsilon q^{\mp ( N + 1 - \varepsilon )} \right)
( X_2 Y_1 + Y_2 X_1 ) + 
\left( 1 - \varepsilon q^{\mp ( N + 1 - \varepsilon)} \right)^2
Y_2 Y_1,
\label{eq:Aexp}
\end{eqnarray}
\begin{eqnarray} 
&&B_{(\pm)} = \hat{R}_{12}^{\pm 1} \hat{R}_{23}^{\pm 1}\nonumber\\ 
&&\phantom{B_{(\pm)}} = {\sf I} - \bigl( 1 + q^{\mp 2}\bigr)
( X_1 + X_2 ) - 
\left( 1 - \varepsilon q^{\mp ( N + 1 - \varepsilon)} \right)
( Y_1 + Y_2 ) + \left( 1 + q^{\mp 2} \right)^2 X_1 X_2\nonumber\\ 
&&\phantom{B_{(\pm)} =} + ( 1 + q^{\mp 2} ) 
\left( 1 - \varepsilon q^{\mp ( N + 1 - \varepsilon )} \right)
( X_1 Y_2 + Y_1 X_2 ) + 
\left( 1 - \varepsilon q^{\mp ( N + 1 - \varepsilon)} \right)^2
Y_1 Y_2. 
\label{eq:Bexp}
\end{eqnarray}  
We also enlist other useful relations:
\begin{equation}
Y_i Y_j Y_i = \left( 1 + \varepsilon [ N - \varepsilon ] 
\right)^{-2} Y_i,
\label{eq:YYY}
\end{equation}
\begin{equation}
Y_i X_j Y_i = \frac{\varepsilon [ N - \varepsilon ] 
\left( [2] + \varepsilon [ N - 1 - \varepsilon ] \right)}
{[2] \left( 1 + \varepsilon [ N - \varepsilon ] \right)^2} Y_i,
\label{eq:YXY}
\end{equation}
where $( i, j ) = ( 1, 2 )$ or $( 2, 1 )$. Employing the above 
constraints and defining the following quantities 
\begin{eqnarray}
&&c = \frac{\varepsilon [ N - \varepsilon + 1 ]}{[2] \left( 1 +
\varepsilon [ N - \varepsilon ] \right)},\qquad d = \left( 1 + 
\varepsilon [ N - \varepsilon ] \right)^{-1},\nonumber\\
&&k = (q - q^{-1})\,\frac{q^{(N + 1 - \varepsilon)/2} - \varepsilon
q^{- (N + 1 - \varepsilon)/2}} {q^{(N - 1 - \varepsilon)/2} + 
\varepsilon q^{- (N - 1 - \varepsilon)/2}}, 
\label{eq:cdkdef}
\end{eqnarray} 
we obtain the promised reduction:
\begin{eqnarray}
&&K_1 = - c \left( J_2 + ( 1 - c ) S_2 \right),\, 
K_2 = ( 1 - c ) \left( J_1 + J_2 + ( 1 - c ) S_2 \right),\, 
K_3 = - c \left( J_1 + ( 1 - c ) S_2 \right),\nonumber\\
&&L_1 = - d \left( J_1 + ( 1 - c ) S_2 \right),\quad 
L_2 = ( 1 - c ) ( 1 - d ) S_2,\quad 
L_3 = - d \left( J_2 + ( 1 - c ) S_2 \right),\nonumber\\ 
&&T_1 = [2]^{-2}\,\left( S_1 + k \left( J_1 + J_2 + ( 1 - c ) S_2 
\right) \right),\qquad\qquad T_2 = d^{2} S_{2}.
\label{eq:tlreduce}
\end{eqnarray}
It is useful to define the ratios of the spectral variables $(w/v)$ 
for the ${\hat R}^{\pm 1}$ matrices as given in (\ref{eq:RRinvbrd}): 
\begin{equation}
f_{\pm} = q^{\pm 1} [2]^{-1} \left( 1 - \varepsilon 
q^{\mp ( N + 1 - \varepsilon)} \right) = - \varepsilon 
q^{\mp ( N + 1 - \varepsilon )} f_{\mp}.
\label{eq:wvrat}
\end{equation}
Then it may be shown that
\begin{eqnarray}
&&c = [2]^{-1} \left(\frac{ q^{-1} f_{+} - q f_{-}}{f_{+} - f_{-}} 
\right),\,\,\, d = \frac{q - q^{-1}}{[2] ( f_{+} - f_{-} )},
\nonumber\\
&&k = - ( q^{2} - q^{-2} )\left(\frac{f_{+} f_{-}}{f_{+} - f_{-}}
\right) = - [2]^{2} d\,f_{+} f_{-}.
\label{eq:cdkf}
\end{eqnarray}
For the special case of the quantum group $SO_q (3)$, where 
$\varepsilon = 1$ and $N = 3$, we obtain the parametric values
\begin{equation}
c = 1 - [2]^{-1},\quad d = ( 1 + [2] )^{-1},\quad 
k = - ( 1 + [2] ) ( 2 - [2] ).
\label{eq:cdkSO3}
\end{equation}   
As a valuable consistency check we carried through the reduction
starting directly with $\varepsilon = 1, N = 3$. Agreement for 
this special case was obtained as the reduction scheme 
(\ref{eq:tlreduce}) was reproduced with the appropriate values of 
the parameters (\ref{eq:cdkSO3}).

\medskip

\noindent{\large{\bf{3.2 Modified braid equation}}}

\medskip

\noindent For the choice of the variables $v = v^{\prime} = 
v^{\prime\prime}$ and $w = w^{\prime} = w^{\prime\prime}$, the 
result (\ref{eq:SORvw}) reduces to
\begin{eqnarray}
&&{\hat R}_{12} (v, w)\,{\hat R}_{23} (v, w)\,{\hat R}_{12} (v, w)
- {\hat R}_{23} (v, w)\,{\hat R}_{12} (v, w)\,{\hat R}_{23} (v, w)
\nonumber\\
&&\quad= v ( 1 + v ) S_1 + w ( 1 + w ) S_2 + v^2 w S_3 + v w^2 S_4
+ v^3 T_1 + w^3 T_2,
\label{eq:mbRvweq}
\end{eqnarray} 
where, using (\ref{eq:tlreduce}) we obtain 
\begin{eqnarray}
&&S_3 = K_1 + K_2 + K_3 = ( 1 -  c ) ( 1 - 3 c ) S_2 +
( 1-2 c ) ( J_1 + J_2 ),\nonumber\\  
&&S_4 = L_1 + L_2 + L_3 = (1 - c ) ( 1 - 3 d ) S_2 - 
d ( J_1 + J_2 ). 
\label{eq:S3S4}
\end{eqnarray} 
Using the reduction (\ref{eq:tlreduce}) for the trilinear 
constructs $T_1$ and $T_2$, we finally obtain the right hand side of 
(\ref{eq:mbRvweq}) as
\begin{equation}
a_1 S_1 + a_2 S_2 + b ( J_1 + J_2 ),
\label{eq:mbrhs} 
\end{equation} 
where the coefficients read
\begin{eqnarray}
&&a_1 = [2]^{-2} v ( v + [2] q ) ( v + [2] q^{- 1} ),\nonumber\\
&&a_2 = w ( 1 + w ) + d^2 w^3 + ( 1 - c ) v^3
\left( [2]^{-2} k + ( 1 - 3 c ) \left( \frac{w}{v} \right) +   
( 1- 3 d ) \left( \frac{w}{v}\right)^2 \right),
\nonumber\\
&&b = - d v^3 \left( \frac{w}{v} - f_{+} \right) 
\left( \frac{w}{v} - f_{-} \right). 
\label{eq:mbcoeff}
\end{eqnarray}
The choice of the `braid values' for the variables $v = - [2] 
q^{\mp 1}, w/v = f_{\pm}$ for the ${\hat R}^{\pm 1}$ matrices 
readily reduces (\ref{eq:mbcoeff}) to $ a_1 = 0,
a_2 = 0, b = 0$. This is obvious for $a_1$ and $b$, whereas the result
for $a_2$ provides a good consistency check.    
 
To obtain the general MBE we now express (\ref{eq:mbrhs}) in terms
of $\left({\hat R} ( v, w )\right)^{\pm 1}$. The projectors now read
\begin{eqnarray}
&&P_{(-)} = \frac{1 + v}{v ( v - w)} \left( \left( {\hat R} ( v, w )
- {\sf I} \right) + ( 1 + w ) \left( {\hat R}^{-1} ( v, w ) - {\sf I}
\right) \right),\nonumber\\
&&P_{(0)} = \frac{1 + w}{w ( w - v)} \left( \left( {\hat R} ( v, w )
- {\sf I} \right) + ( 1 + v ) \left( {\hat R}^{-1} ( v, w ) - {\sf I}
\right) \right).
\label{eq:PMPN}
\end{eqnarray}
In passing we mention that by implementing the `braid values'
(\ref{eq:brdval}) of the variables $( v, w )$ the projectors
$P_{(+)}$ and $P_{(0)}$ may also be expressed in terms of the braid
matrices ${\hat R}^{\pm 1}$. Now we can express the constructs 
$S_1, S_2, J_1$ and $J_2$ defined in (\ref{eq:SOdef}) in terms of 
the generalized braid matrices ${\hat R}_{12}^{\pm 1} ( v, w )$ 
and ${\hat R}_{23}^{\pm 1} ( v, w )$. Substituting these results 
in (\ref{eq:mbrhs}) we obtain the general MBE stated in 
(\ref{eq:BCDmbe}), where, with the values of $( a_1, a_2, b )$ 
given in (\ref{eq:mbcoeff}), the coefficients read
\begin{eqnarray}
&&c_1 = \frac{( 1 + v ) ( 1 + w )}{v w ( v - w )} \left( \frac{w}
{1+ w} a_1 - \frac{v}{1 + v} a_2 + 2 b \right),\nonumber\\
&&c_2 = \frac{( 1 + v ) ( 1 + w )}{v w ( v - w )} 
( w a_1 - v a_2 - 2 b ),\qquad
c_3 = - \frac{( 1 + v ) ( 1 + w )}{v w ( v - w )} b.
\label{eq:c123}
\end{eqnarray}
For the `braid values' (\ref{eq:brdval}) of the variables $v$ and 
$w$ the coefficients $( c_1, c_2, c_3 )$ vanish. 

We now discuss the following special cases depending on values
of the parameters $v$ and $w$.
\begin{itemize}

\item  For {\em arbitrary} non-zero $v$, and for $w = f_{\pm} v$, 
we, {\em via} (\ref{eq:mbcoeff}) and (\ref{eq:c123}), obtain 
$c_3 = 0$. Hence the right hand side of the MBE (\ref{eq:BCDmbe}) 
reduces to the first two terms. Using the characteristic equation
(\ref{eq:BCDfact}) the right hand side of the MBE may now be 
expressed in terms of the generalized braiding operators
$\left( {\hat R}_{12} ( v, w ) - {\hat R}_{23} ( v, w ) \right)$ and
$\left( {\hat R}_{12}^2 ( v, w ) - {\hat R}_{23}^2 ( v, w ) \right)$. 

\item  For $v = 0$ and an arbitrary value of $w$ we have an 
interesting case deserving a detailed treatment that is provided
in Sec. 4. Here we just note that, for these parametric values,
the operator ${\hat R} ( 0, w )$ satisfies a {\em quadratic} rather
than a cubic characteristic equation; and the right hand side of the 
MBE (\ref{eq:BCDmbe}) is now proportional to    
$\left( {\hat R}_{12} ( 0, w ) - {\hat R}_{23} ( 0, w ) \right)$. 

\end{itemize}

\medskip

\noindent{\large{\bf 3.3 Baxterization}} 

\medskip

\noindent As was noted in Sec.1, the MBE and the Baxterization are 
two complementary aspects of the generalized braid operator 
${\hat R} ( v, w )$ with its spectral decomposition given in 
(\ref{eq:BCDRvw}). Having formulated the MBE, we now turn to 
Baxterization of the braid operator. 

{}From (\ref{eq:SORvw}) and (\ref{eq:tlreduce}) we obtain
\begin{eqnarray}
&&{\hat R}_{12} ( v, w ) {\hat R}_{23} ( v^{\prime}, w^{\prime} )
{\hat R}_{12} ( v^{\prime\prime}, w^{\prime\prime} ) -
{\hat R}_{23} ( v^{\prime\prime}, w^{\prime\prime} ) 
{\hat R}_{12} ( v^{\prime}, w^{\prime} ) 
{\hat R}_{23} ( v, w )\nonumber\\ 
&&\qquad = {\sf a}_1 S_1 + {\sf a}_2 S_2 + {\sf b}_1 J_1 + 
{\sf b}_2 J_2,
\label{eq:SSJJ}
\end{eqnarray}
where, with the parameters $c, d$ and $k$ given by 
(\ref{eq:cdkdef}), the coefficients read 
\begin{eqnarray}
&&{\sf a}_1 = \bigl( v + v^{\prime\prime} + v v^{\prime\prime} 
\bigr) - \bigl( 1 - [2]^{-2} v v^{\prime\prime} \bigr) 
v^{\prime},\nonumber\\
&&{\sf a}_2 = \bigl( w + w^{\prime\prime} + w w^{\prime\prime} 
\bigr) - \bigl( 1 - d^2\,w w^{\prime\prime} \bigr) w^{\prime} +
( 1 - c )\,\bigl( ( 1 - c )\,v v^{\prime\prime} w^{\prime} -
c\,( v w^{\prime\prime} + v^{\prime\prime} w ) v^{\prime} 
\bigr) \nonumber\\
&&\phantom{{\sf a}_2 = } + ( 1 - c ) \bigl( ( 1 - d )\, 
w w^{\prime\prime} v^{\prime} 
- d\,( v w^{\prime\prime} + v^{\prime\prime} w ) w^{\prime} \bigr) 
+ ( 1 - c ) k [2]^{-2} v v^{\prime\prime} v^{\prime},\nonumber\\
&&{\sf b}_1 = ( v w^{\prime} - w v^{\prime} ) 
- c\,w v^{\prime\prime} v^{\prime} 
+ ( 1 - c )\,v v^{\prime\prime} w^{\prime}
- d\,w v^{\prime\prime} w^{\prime} 
+ k [2]^{-2} v v^{\prime\prime} v^{\prime},\nonumber\\
&&{\sf b}_2 = ( v^{\prime\prime} w^{\prime} 
- w^{\prime\prime} v^{\prime} )
- c\,v w^{\prime\prime} v^{\prime} 
+ ( 1 - c )\,v v^{\prime\prime} w^{\prime} 
- d\,v w^{\prime\prime} w^{\prime} 
+ k [2]^{-2} v v^{\prime\prime} v^{\prime}.  
\label{eq:aabb}
\end{eqnarray}
For Baxterization the variables $( v, v^{\prime\prime}, 
v^{\prime} )$ and $( w, w^{\prime\prime}, w^{\prime} )$ are 
reexpressed as $( v (x), v (y), v (xy) )$ and $ ( w (x), w (y), 
w (xy) )$ respectively; and then the functional equations
\begin{equation}
{\sf a}_i ( x, y ) = 0,\qquad {\sf b}_i ( x, y ) = 0,
\qquad i = ( 1, 2 )
\label{eq:abxy}
\end{equation}
obtained {\em via} ({\ref{eq:aabb}) are solved. 

We start by noting a crucial constraint on the functional solutions 
to be implemented systematically. In (\ref{eq:SSJJ}) the 
${\hat R} ( v, w ) $ matrix may be replaced with its inverse  
${\hat R}^{-1} ( v, w ) $ given in (\ref{eq:invRvw}). Hence for 
consistency our functional solutions must be such that 
\begin{equation}
- \frac{v (x)}{1 + v (x)} = v (x^{\prime}),\qquad 
- \frac{w (x)}{1 + w (x)} = w (x^{\prime}),
\label{eq:xprime}
\end{equation}
where $x^{\prime}$ is some suitable function of $x$. We have 
chosen the parametrization 
\begin{equation}
x^{\prime} = x^{-1}
\label{eq:xrel}
\end{equation} 
such that the validity of the constraint
\begin{equation}
{\hat R} (x) {\hat R} (x^{-1}) = {\sf I}
\label{eq:RR}
\end{equation}
is assured. The functions $f(x)$ and $g(x)$ being suitably 
well-defined, the relations (\ref{eq:xprime}) and (\ref{eq:xrel}) 
now imply
\begin{equation}
v (x) = \frac{f ( x )}{f ( x^{-1} )} - 1,\qquad
w (x) = \frac{g ( x )}{g ( x^{-1} )} - 1.
\label{eq:vwx}
\end{equation}
It will, in fact, thus be necessary to solve only {\em two} 
functional equations {}from the set (\ref{eq:abxy}):
\begin{equation} 
{\sf a}_1 ( x, y ) = 0,\qquad {\sf b}_1 ( x, y ) = 0.
\label{eq:a1b1xy}
\end{equation}
As for the function ${\sf b}_2$ we note that the following 
exchange relation holds:
\begin{equation}
{\sf b}_2 ( x, y ) = {\sf b}_1 ( y, x ).
\label{eq:b1b2}   
\end{equation}  
Moreover the solutions of the functional equations 
(\ref{eq:a1b1xy}), along with (\ref{eq:vwx}), will be seen to fix 
the functions $v (x)$ and $w (x)$ {\em completely with all the 
parameters determined}. Then the much more complicated functional 
equation for ${\sf a}_2 ( x, y )$ does not have to be `solved' at 
all. The consistency of the scheme may be verified by checking 
that ${\sf a}_2 ( x, y )$ indeed vanishes by implementing $v (x)$ 
and $w (x)$ already available. So now we consider the equations 
in (\ref{eq:a1b1xy}). The first equation in (\ref{eq:a1b1xy}) 
implies (\ref{eq:vbaxrel}). It will turn out to be quite useful to 
solve a more general function equation with arbitrary $\lambda$ 
given by
\begin{equation}
u ( x y ) = \frac{u (x) + u (y) + u (x) u (y)}{1 - \lambda^{2}
u (x) u (y)}.
\label{eq:uxyfun}
\end{equation}    
By direct substitution it may be verified that the general solution 
maintaining the structure (\ref{eq:vwx}), namely
\begin{equation}
u (x) = \frac{\alpha ( x )}{\alpha ( x^{-1} )} - 1
\label{eq:uxsol}
\end{equation}
is given by
\begin{equation}
\alpha (x) = - ( x - x^{-1} ) \pm \sqrt{1 - 4 \lambda^{2}}
( x + x^{-1} ).
\label{eq:alphax}
\end{equation} 
For all cases of interest to us the square root will turn out to 
be real. The lower sign before the root in the right hand side of 
(\ref{eq:alphax}) corresponds to $u ( x^{-1} )$    with the upper 
sign. So it is sufficient to consider only one, say, the upper sign 
before the root, since, more generally, replacing of $x$ by $x^p$ 
in (\ref{eq:vwx}) does not change the essential functional form. 
Now setting $\lambda = [2]^{-1}$, we obtain the solution of the 
first equation in (\ref{eq:vwx}) with the function $f (x)$ given by
\begin{equation}
f (x) = - ( x - x^{-1} ) + ( q - q^{-1} ) [2]^{-1} ( x + x^{-1} ).
\label{eq:fxsol}
\end{equation}
Apart from the basic structures (\ref{eq:uxsol}) and 
(\ref{eq:alphax}), the following form is also of interest:
\begin{equation}
u (x) = \beta \left(\frac{x^2 - 1}{x^2 + \beta + 1} \right),\qquad
\hbox{where}\,\,\beta = - 2\left( 1 + \sqrt{1 - 4 \lambda^2} 
\right)^{-1}.
\label{eq:altux}    
\end{equation}   
    
Substituting $x = \hbox{exp} (\theta) , \hbox{tanh}\,\gamma = 
\sqrt{1 - 4 \lambda^{2}}$ we obtain the elegant form of additive 
Baxterization, as seen in the context of (\ref{eq:adbax}):
\begin{equation}
u (x) \sim u ( \theta ) = \frac{\hbox{sinh} ( \gamma - \theta )}
{\hbox{sinh} ( \gamma + \theta )} - 1.
\label{eq:uxthta}
\end{equation}
For the function $v (x) \sim v ( \theta )$, we obtain, after 
setting $q = \hbox{exp} (h)$, the solution (\ref{eq:vthta}). 

Now to solve the second equation in (\ref{eq:a1b1xy}) we 
proceed in successive steps as follows. Expressing the parameters
$( c, d, k )$ in terms of the ratios $f_{\pm}$ as given in 
(\ref{eq:cdkf}), we define 
\begin{equation}
{\sf b}_1 \equiv - \frac{v v^{\prime} v^{\prime\prime}}
{[2] ( f_{+} - f_{-} )} F,
\label{eq:b1F} 
\end{equation}
where
\begin{eqnarray}
&&F = ( q - q^{-1} ) \frac{w}{v} \frac{w^{\prime}}{v^{\prime}} +
[2] \frac{f_{+} - f_{-}}{v^{\prime\prime}} \left(\frac{w}{v} - 
\frac{w^{\prime}}{v^{\prime}} \right) + ( q^{-1} f_{+} - q f_{-} ) 
\frac{w}{v}\nonumber\\ 
&&\phantom{F =} - ( q f_{+} - q^{-1} f_{-} ) 
\frac{w^{\prime}}{v^{\prime}} + ( q - q^{-1} ) f_{+} f_{-}. 
\label{eq:Fexp}
\end{eqnarray}
In (\ref{eq:Fexp}) we have assumed that the spectral function
$v \neq 0$. Now defining 
\begin{equation}
U = [2] d ( w v^{-1} - f_{+} ),
\label{eq:Udef}
\end{equation}
we obtain after simplification
\begin{equation}
F = ( q - q^{-1} ) \frac{U U^{\prime}}{( [2] d)^2} 
\left( 1 + \left( q + \frac{[2]}{v^{\prime\prime}} \right) 
\frac{1}{U^{\prime}} - \left( q^{-1} + \frac{[2]} 
{v^{\prime\prime}} \right) \frac{1}{U} \right).
\label{eq:FUdef}
\end{equation}          
In the context of (\ref{eq:Udef}) we note that by translating 
the ratio $( w/v )$ by $f_{-}$ the same functional solutions 
are finally obtained. Substituting for the function $v$ the form 
corresponding to (\ref{eq:altux}), {\em i.e.} 
\begin{equation} 
v^{\prime\prime} = v ( y ) = - [2] q^{-1} ( y^2 -1 ) 
( y^2 - q^{-2} )^{-1},
\label{eq:vyval}
\end{equation}
we find that the requirement for $F = 0$ reads
\begin{equation}
U^{-1} ( x y ) + ( q - q^{-1} )^{-1} = \left( U^{-1} ( x ) +
( q - q^{-1} )^{-1} \right) y^2. 
\label{eq:Ufneq}
\end{equation}
Hence, with an as yet arbitrary parameter $\delta$ that is to be 
fixed immediately afterwards, we obtain the solution as 
\begin{equation}
U^{-1} ( x ) + ( q - q^{-1} )^{-1} = 
\delta ( q - q^{-1} )^{-1} x^2.
\label{eq:Usol}
\end{equation}
Now the spectral function $w ( x )$ may be solved in terms of the 
already known solution $v ( x )$: 
\begin{equation}
w ( x ) = ( f_{+} \delta x^2 - f_{-} ) 
( \delta x^2 - 1 )^{-1}\,v ( x ). 
\label{eq:wxsol}
\end{equation}
Using the crucial constraint (\ref{eq:vwx}) on $w ( x )$  now it 
may be shown that there are {\em only two permitted values} of the
parameter $\delta$, which, in turn, generate the corresponding 
solutions of the spectral function $w ( x )$:    
\begin{equation}
\delta = - \varepsilon q^{( N + 1 - \varepsilon )}
\quad\Rightarrow\quad w ( x ) = \frac{x + \varepsilon 
q^{( N + 1 - \varepsilon )} x^{-1}} {x^{- 1} + \varepsilon 
q^{( N + 1 - \varepsilon )} x} - 1,
\label{eq:wx1}
\end{equation}
\begin{equation}
\delta = \varepsilon q^{( N - 1 - \varepsilon )}\quad
\Rightarrow\quad w ( x ) = \frac{\left( q x - \varepsilon 
q^{( N - \varepsilon )} x^{-1}\right) ( q^{-1} x - q x^{-1})}  
{\left( q x^{-1} - \varepsilon q^{( N - \varepsilon )} x \right) 
( q^{-1} x^{-1} - q x )} - 1.  
\label{eq:wx2}
\end{equation} 
At this stage our procedure of Baxterization of the braid operator 
is complete. In view of the discussion following (\ref{eq:b1b2}) 
we need not solve the functional equations corresponding to 
${\sf a}_2 ( x, y )$ and ${\sf b}_2 ( x, y )$. The function 
$w ( x )$ corresponding to the additive form of Baxterization may 
be obtained by setting $x = \hbox{exp} ( \theta ), 
q = \hbox{exp} ( h )$ in the solutions (\ref{eq:wx1}) and 
(\ref{eq:wx2}). For the quantum group $SO_q ( N )$, with 
$\varepsilon = 1$, these solutions assume the form (\ref{eq:wSO1}) 
and (\ref{eq:wSO2}) respectively. Similarly for the quantized 
symplectic group $Sp_q ( N )$, where $\varepsilon = - 1$ and 
$N = 2 n$, they are given by (\ref{eq:wSp1}) and (\ref{eq:wSp2}) 
respectively.     

Finally, expressing the projection operators $P_{(-)}$ and 
$P_{(0)}$ in terms of the braid matrices ${\hat R}^{\pm 1}$, as 
discussed following equation (\ref{eq:PMPN}), we obtain
\begin{equation}
{\hat R} ( x ) = l_{(0)} ( x ) {\sf I} + l_{(+)} ( x ) {\hat R}
+ l_{(-)} ( x ) {\hat R}^{-1}.
\label{eq:RL}
\end{equation}
Denoting the `braid values' (\ref{eq:brdval}) of the spectral 
variables corresponding to the ${\hat R}$ matrix by $( v_0, w_0 )$ 
respectively, the coefficients in (\ref{eq:RL}) may be listed as
\begin{eqnarray}
&&l_{(0)} ( x ) = 1 - ( 1 + v_0 ) ( 2 + w_0  ) v_{0}^{-1} 
( v_0 - w_0 )^{-1} \,v ( x )\nonumber\\
&&\phantom{l_{(0)} ( x ) =} - ( 2 + v_0 ) ( 1 + w_0 ) w_{0}^{-1} 
( w_0 - v_0 )^{-1} \,w ( x ),\nonumber\\
&&l_{(+)} ( x ) = ( 1 + v_0 ) v_{0}^{-1} ( v_0 - w_0 )^{-1}\, 
v ( x ) + ( 1 + w_0 ) w_{0}^{-1} ( w_0 - v_0 )^{-1}\, 
w ( x ),\nonumber\\
&&l_{(-)} ( x ) = ( 1 + v_0 ) ( 1 + w_0 ) ( v_{0}^{-1} 
( v_0 - w_0 )^{-1} \, v ( x ) + w_{0}^{-1} 
( w_0 - v_0 )^{-1} w ( x ) ).
\label{eq:Lxrel}
\end{eqnarray} 
Consistent with our normalization, the relation  $l_{(0)} ( x ) + 
l_{(+)} ( x ) + l_{(-)} ( x ) = 1$ is maintained.  

\sect{A new class of solutions of the braid equation}
So far we have been studying the general case $( v \neq 0, 
w \neq 0 )$ of the spectral variables. We note that for the choice 
$w = 0$, there is no solution of the braid equation with nontrivial 
$v$ as more than one functional constraints on the single spectral 
parameter, or the function $v( x )$ in the Baxterized case, need to 
be satisfied.

For the alternate choice $v = 0$, it is, however, {\em possible} 
to obtain a remarkable class of generalized braid matrices 
${\hat R} ( w )$ depending on the single variable $w$. This was 
already studied in the context of the quantum group $SO_q (3)$ 
in \cite{C02}. Here we present this class for general quantum 
groups $SO_q ( N )$ and $Sp_q ( N )$. Denoting the generalized 
braid matrix ${\hat R} ( 0, w )$ by ${\hat R} ( w )$, we write
\begin{equation}
{\hat R} ( w ) = {\sf I} + w P_{(0)},\qquad {\hat R}^{-1} ( w ) = 
{\sf I} - w ( 1 + w )^{-1} P_{(0)},
\label{eq:RRinvw}
\end{equation}
where the matrix ${\hat R} ( w )$ satisfies a quadratic 
characteristic equation
\begin{equation}
( {\hat R} ( w ) - {\sf I} )\,( {\hat R} ( w ) - 
( 1 + w )\,{\sf I} ) = 0.
\label{eq:Rwchar}
\end{equation}
Instead of restricting our preceding general formulae to the 
special case $v = 0$ it is much simpler to start again with 
(\ref{eq:SORvw}) which, in conjunction with the last relation in 
(\ref{eq:tlreduce}), now gives 
\begin{eqnarray}
&&{\hat R}_{12} ( w ) {\hat R}_{23} ( w ) {\hat R}_{12} ( w ) -
{\hat R}_{23} ( w ) {\hat R}_{12} ( w ) {\hat R}_{23} ( w )
\nonumber\\      
&&\,\, = ( w ( 1 + w ) + d^2 w^3 ) S_2 = ( 1 + w + d^2 w^2 )
\left( {\hat R}_{12} ( w ) - {\hat R}_{23} ( w ) \right).
\label{eq:mbew}
\end{eqnarray}      
Thus we have already obtained the MBE with a simple linear 
structure in the generalized braid matrix on the right hand side. 
But as the coefficient on the right of (\ref{eq:mbew}) may be 
factorized as
\begin{equation}
d^2 w^2 + w + 1 = d^2 ( w - w_{+} ) ( w - w_{-} ),\qquad
\,\,w_{\pm} = 2^{-1} d^{-2} \left( - 1 \pm \sqrt{ 1 - 4 d^2 } \right),
\label{eq:wroot}
\end{equation}      
we obtain two `non-modified' {\em new braid matrices}:
\begin{equation}       
{\hat R}_{12} ( w_{\pm} ) {\hat R}_{23} ( w_{\pm} ) 
{\hat R}_{12} ( w_{\pm} ) - {\hat R}_{23} ( w_{\pm} ) 
{\hat R}_{12} ( w_{\pm} ) {\hat R}_{23} ( w_{\pm} ) = 0.
\label{eq:newbr}
\end{equation}
But the constraint
\begin{equation}
{\hat R} ( w_{+} ) {\hat R} ( w_{-} ) = {\sf I},
\label{eq:oneR}
\end{equation}
observed using (\ref{eq:RRinvw}) and (\ref{eq:wroot}) indicate 
the presence of only {\em one} braid matrix and its inverse.
Following the expression of $d$ in (\ref{eq:cdkdef}) we  also 
note that the quantity under the radical sign in (\ref{eq:wroot}) 
satisfies $1 > ( 1 - 4 d^2 ) > 0$ for all values of $( N, 
\varepsilon )$ relevant to us. Hence the square root is real. 

The prescription of Baxterization of this class of 
${\hat R} ( w )$ matrices is immediately obtained {}from 
(\ref{eq:uxyfun})-(\ref{eq:alphax}). Only relevant functional 
equation ${\sf a}_2 ( x, y ) = 0$ in the set 
(\ref{eq:abxy}) now yields
\begin{equation}
w ( x y ) = \frac{w ( x ) + w ( y ) + w ( x ) w ( y )}{1 - d^2 
w ( x )\,w ( y )},
\label{eq:newbrfun}
\end{equation}
whereas the other functional equations in the set (\ref{eq:abxy})
trivially vanishes for the case $v = 0$. The solution of the equation
(\ref{eq:newbrfun}) read 
\begin{equation}
w ( x ) = \frac{- ( x - x^{-1} ) + \sqrt{1 - 4 d^2} ( x + x^{-1} )} 
{( x - x^{-1} ) + \sqrt{1 - 4 d^2} ( x + x^{-1} )} - 1.
\label{eq:newbrwx}
\end{equation} 
Setting $x = \hbox{exp} ( \theta )$ and $\hbox{tanh}\,\eta = 
\sqrt{1 - 4 d^2}$, we, following (\ref{eq:uxthta}), obtain
\begin{equation}
w ( x ) \sim w ( \theta ) = \frac{\hbox{sinh} ( \eta - \theta )} 
{\hbox{sinh} ( \eta + \theta )} - 1.
\label{eq:newbrwth}
\end{equation} 

For the quantum group $SO_q (3)$ the generalized braid operator 
${\hat R} ( w )$ in the present class was explicitly presented 
before in \cite {C02}. Here we present the general prescription
regarding the quantum groups $SO_q ( N )$ and $Sp_q ( N )$. To this
end we include the standard and well-known construction of the 
projectors for the sake of completeness. Let $( \rho_1, \rho_2,
\cdots,\rho_N )$ denote the $N$-tuple with the following assignments 
for the respective quantum groups:
\begin{eqnarray}
&&\left( n - \frac{1}{2}, n - \frac{3}{2},\cdots,\frac{1}{2}, 0,
- \frac{1}{2},\cdots, - n + \frac{1}{2} \right)\quad \hbox{for}\,
SO_q ( 2 n + 1 ),\nonumber\\
&&( n - 1, n - 2,\cdots, 1, 0, 0, -1,\cdots, - n + 1 )\quad 
\hbox{for}\,SO_q ( 2 n ),\nonumber\\
&&( n, n - 1,\cdots,1, -1,\cdots, - n )\quad \hbox{for}\,
Sp_q ( 2 n ).
\label{eq:rhoval}
\end{eqnarray}
We also fix the values $\{ \epsilon_i = 1 | i = 1,\cdots, N \}$ 
for the quantum group $SO_q ( N )$; and $\{ \epsilon_i = 1,\,\,
\hbox{if}\,\,i \leq n, \epsilon_i = - 1,\,\,\hbox{if}\,\,i > n \}$ 
for the quantum group $Sp_q ( 2 n )$. As usual $E_{ij}$ denotes 
the $N \times N$ matrix where the only non-zero element 
$( i, j )$ equals unity. With these notations the projectors 
are given by 
\begin{equation}
P_{(0)} = ( 1 + [ N - 1 ] )^{-1} \sum^{N}_{i, j = 1} q^{( \rho_i
- \rho_j )}\,\, E_{i^{\prime}, j} \otimes E_{i, j^{\prime}}
\label{eq:PSO}
\end{equation}
and
\begin{equation}
P_{(0)} = ( 1 - [ N + 1 ] )^{-1} \sum^{N}_{i, j = 1} q^{( \rho_i
- \rho_j)}\,\,\epsilon_i \epsilon_j\,\, E_{i^{\prime}, j} \otimes 
E_{i, j^{\prime}}
\label{eq:PSp}
\end{equation}
for the quantum groups $SO_q(N)$ and $Sp_q (N),\,( \hbox{where}\,
N = 2 n )$ respectively. In the above equations we have used the
notation $i^{\prime} = N + 1 - i$.  
We have recapitulated the above standard prescription for the 
additional purpose of displaying the braid operators in the 
$q = 1$ limit.  

We enlist below the lowest dimensional cases as examples. The 
projector $P_{(0)}$ in the $9 \times 9$ vector representation 
for the quantum group $SO_q (3)$ is given by
\begin{eqnarray}
&&( 1 + [2] ) P_{(0)} = q^{-1} E_{11} \otimes E_{33} +
q^{- \half} E_{12} \otimes E_{32} + E_{13} \otimes E_{31} +
q^{- \half} E_{21} \otimes E_{23} + E_{22} \otimes E_{22} 
\nonumber\\
&&\phantom{( 1 + [2]) P_{(0)} =}  
+ q^{\half} E_{23} \otimes E_{21} + E_{31} \otimes E_{13} 
+ q^{\half} E_{32} \otimes E_{12}
+ q E_{33} \otimes E_{11},
\label{eq:SO3P}
\end{eqnarray}
whereas $16 \times 16$ vector representations for the quantum 
groups $SO_q (4)$ and $Sp_q (4)$ read as
\begin{eqnarray}
&&( 1 + [3] ) P_{(0)} = q^{-2} E_{11} \otimes E_{44} 
+ q^{- 1} E_{12} \otimes E_{43} + q^{-1} E_{13} \otimes E_{42} 
+ E_{14} \otimes E_{41}\nonumber\\
&&\phantom{( 1+ [3] ) P_{(0)} =} + q^{-1}E_{21} \otimes 
E_{34} + E_{22} \otimes E_{33} + E_{23} \otimes E_{32} 
+ q E_{24} \otimes E_{31} \nonumber\\ 
&&\phantom{( 1 + [3] ) P_{(0)} =} + q^{-1} E_{31} \otimes 
E_{24} + E_{32} \otimes E_{23} + E_{33} \otimes E_{22} 
+ q E_{34} \otimes E_{21}\nonumber\\ 
&&\phantom{( 1 + [3] ) P_{(0)} =} + E_{41} \otimes E_{14} 
+ q E_{42} \otimes E_{13} + q E_{43} \otimes E_{12} 
+ q^2 E_{44} \otimes E_{11} 
\label{eq:SO4P}
\end{eqnarray}
and
\begin{eqnarray}
&&( [5] - 1 ) P_{(0)} = 
q^{-4} E_{11} \otimes E_{44} + q^{-3} E_{12} \otimes E_{43} 
- q^{-1} E_{13} \otimes E_{42} - E_{14} \otimes E_{41}\nonumber\\ 
&&\phantom{( [5] -1 ) P_{(0)} =} 
+ q^{-3} E_{21} \otimes E_{34} + q^{-2} E_{22} \otimes E_{33} 
- E_{23} \otimes E_{32} - q E_{24} \otimes E_{31}\nonumber\\  
&&\phantom{( [5] -1 ) P_{(0)} =} 
- q^{-1} E_{31} \otimes E_{24} - E_{32} \otimes E_{23} 
+ q^2 E_{33} \otimes E_{22} + q^3 E_{34} \otimes E_{21}\nonumber\\
&&\phantom{( [5] - 1 ) P_{(0)} =} 
- E_{41} \otimes E_{14} - q E_{42} \otimes E_{13} 
+ q^3 E_{43} \otimes E_{12} + q^4 E_{44} \otimes E_{11}
\label{eq:Sp4P}
\end{eqnarray} 
respectively.

{\em For this class, directly setting $q = 1$ we still get 
nontrivial braid matrices, the MBE satisfied by them, and also 
the Baxterized forms of these braid matrices}. ( We note that we 
are not considering the so-called quasi-classical limit obtained 
as the coefficients of the terms linear in $h$, while 
implementing series expansion in $h\,\equiv\,\hbox{ln}\,q$. We 
simply set $q = 1$ in the relevant quantities, and thus obtain 
the terms of zero-order in $h$ ). The corresponding example for 
the quantum group $SO_q (3)$ was discussed in \cite{C02}. 

In the classical $q = 1$ limit, we denote the projectors by 
${\hat P}$ to avoid confusion. Using (\ref{eq:PSO}) we now 
obtain the projector for the classical group $SO ( N )$ as 
\begin{equation}
{\hat P}_{(0)} = N^{-1} \sum_{i, j = 1}^{N} E_{i^{\prime}, j} 
\otimes E_{i, j{\prime}},
\label{eq:clP0SO}
\end{equation}
while the projector for the classical symplectic group 
$Sp (N )$ reads 
\begin{equation}
{\hat P}_{(0)} = - N^{-1} \sum_{i, j = 1}^{N} \epsilon_i 
\epsilon_j E_{i^{\prime}, j} \otimes E_{i, j^{\prime}}.
\label{eq:clP0Sp}
\end{equation}
Moreover the parameters, {\em via} (\ref{eq:cdkdef}), 
(\ref{eq:wroot}) and (\ref{eq:newbrwth}), are obtained 
in the limit $ q = 1$ as
\begin{eqnarray}
&&d\,|_{( q = 1 )} \equiv {\hat d} = \varepsilon\,N^{-1},\qquad
\eta\,|_{( q = 1 )} \equiv {\hat {\eta}} = \hbox{arctanh}\,
( N^{-1} \sqrt{ N^2 -4 } ),\nonumber\\ 
&&w_{( \pm )}|_{( q = 1 )} \equiv {\hat w}_{( \pm )} = \frac{N}{2} 
\left( - N \pm \sqrt{ N^2 - 4} \right). 
\label{eq:q1dweta}
\end{eqnarray}
Using the above value of ${\hat d}$ in (\ref{eq:mbew}) we may 
directly obtain the MBE in the $q = 1$ case. In the Baxterized 
function (\ref{eq:newbrwth}) we now use the value of 
${\hat {\eta}}$ given in (\ref{eq:q1dweta}), and thereby obtain 
the limiting structure
\begin{equation}
w ( \theta )\,|_{( q = 1 )} = {\hat w} ( \theta ) = \frac
{\hbox{sinh}\,( {\hat {\eta}} -\theta )}{\hbox{sinh} 
( {\hat {\eta}} + \theta ) } - 1.
\label{eq:q1wthta}
\end{equation}
We now focus on the $q = 1$ limiting case of the braid solution 
described in (\ref{eq:newbr}). The braid matrices in this limit 
has the structure 
\begin{equation}
{\hat R} ( {\hat w}_{( \pm )} ) = {\sf I} + {\hat w}_{( \pm )} 
{\hat P}_{(0)}. 
\label{eq:q1bm}
\end{equation}   
Even in the said $q = 1$ limit, these matrices, on account of the 
characteristic equation (\ref{eq:Rwchar}), satisfy a nontrivial  
Hecke condition. These braid matrices are not of co-boundary 
type, which may be obtained by twisting the identity operator. 
The situation may be profitably contrasted with the $q = 1$
limiting behavior of the general case $( v \neq 0 )$ of braid 
matrices satisfying (\ref{eq:be}). Setting the limit $q =1$  
in the `braid values' (\ref{eq:brdval}) of the spectral 
variables we obtain  
\begin{equation}
v\,|_{( q = 1 )} \equiv {\hat v} = -2,\qquad 
w\,|_{( q = 1 )} \equiv {\hat w} = - ( 1 - \varepsilon ).
\label{eq:vwhat}
\end{equation}
Incidentally, identical conditions in the said classical limit
$h = 0$ may be derived {}from the Baxterized version of the 
spectral variables given in (\ref{eq:vthta})-(\ref{eq:wSp2}).
For the values of the spectral variables given in (\ref{eq:vwhat})
the braid matrix for the quantum group $SO_q ( N )$ assumes the form
\begin{equation}
{\hat R} ( {\hat v}, {\hat w} ) = {\sf I} - 2 P_{(-)}, 
\label{eq:RclSO}
\end{equation}
whereas the the braid matrix for the quantum group $Sp_q ( N )$
reads       
\begin{equation}
{\hat R} ( {\hat v}, {\hat w} ) = {\sf I} - 2 P_{(-)} - 2 P_{(0)} 
= 2 P_{(+)} - {\sf I}.
\label{eq:RclSp}
\end{equation}
The above two braid matrices satisfy the condition 
$\left( {\hat R} ( {\hat v}, {\hat w} ) \right)^2 = {\sf I}$, 
typical of twisted identity matrices. Finally we note that here 
the result parallel to (\ref{eq:Rx}) is 
\begin{equation}
{\hat R} ( w ) = \frac{ ( w - w_{-} ) {\hat R} ( w_{+} ) - 
( w - w_{+} ) {\hat R} ( w_{-} )}{( w_{+} - w_{-} )},
\label{eq:Rwpwm}
\end{equation}
where $w_{\pm}$, {\em via} (\ref{eq:wroot}), is given by $w_{\pm} =
- \left( 1 + \hbox{exp} ( \mp 2 \eta ) \right)$.      

\sect{Diagonalization of the modified braid matrices and their
corresponding eigenvectors}

We present below explicitly the $4 \times 4, 9 \times 9$ and
$16 \times 16$ matrices which, through similarity transformations,
diagonalize respectively the matrices ${\hat R} ( v )$ for the 
quantum group $GL_q (2)$, and ${\hat R} ( v, w )$ for the quantum 
groups $SO_q (3)$ and $SO_q (4)$. Explicit constructions for 
arbitrary dimensions $N \times N$ is beyond the scope of the 
present work. But we start by considering, for deeper 
understanding,  certain aspect of the problem for the general 
quantized orthogonal group $SO_q ( N )$. The case for the general 
linear quantum group $GL_q ( N )$ is, as usual, much simpler. In 
the diagonalization process, as throughout the present work, the 
projectors involved in the spectral decomposition of the 
generalized braid matrices play essential roles. Interest in the 
results obtained will be discussed at the end. 

Traces of the projectors satisfying the completeness property 
\begin{equation}
P_{( + )} + P_{( - )} + P_{( 0 )} = {\sf I}_{N^2 \times N^2},   
\label{eq:TrP}
\end{equation}
are given by
\begin{equation}
Tr\,P_{(+)} = \half N ( N + 1 ) - \half ( \varepsilon + 1 ),\,\,
Tr\,P_{(-)} = \half N ( N - 1 ) + \half ( \varepsilon - 1 ),\,\,
Tr\,P_{(0)} = 1.   
\label{eq:Trexpl}
\end{equation}
The sum of the above three traces, consistent with (\ref{eq:TrP}), 
is $N^2$. We consider below only the examples of quantum orthogonal
algebras corresponding to the choice $\varepsilon = 1$. Moreover, 
as in the previous sections we will only consider the projectors 
$P_{(-)}$ and $P_{(0)}$ in view of (\ref{eq:TrP}). 

A projector when diagonalized can have for each diagonal element
{\em either zero or unity}. The number of unit elements are fixed 
by trace. The diagonal elements can be permuted by successive 
similarity transformations. So we can choose a suitable canonical 
ordering for them as follows. Denoting a diagonalized projector 
$P$ by $P^{({\sf d})}$, we choose the ordering presented below. 
For the quantum group $SO_q (3)$ the diagonalized projectors read
\begin{eqnarray}
&&P_{(0)}^{({\sf d})} = ( 1, 0, 0, 0, 0, 0, 0, 0, 
0 )_{({\sf diagonal})},\nonumber\\
&&P_{(-)}^{({\sf d})} = ( 0, 1, 1, 1, 0, 0, 0, 0, 
0 )_{({\sf diagonal})},
\label{eq:PdgSO3}
\end{eqnarray}
whereas these diagonalized projectors for the quantum group 
$SO_q (4)$ are given by
\begin{eqnarray}
&&P_{(0)}^{({\sf d})} = ( 1, 0, 0, 0, 0, 0, 0, 0, 0, 0, 0, 0, 0, 0, 
0, 0 )_{({\sf diagonal})},\nonumber\\
&&P_{(-)}^{({\sf d})} = ( 0, 1, 1, 1, 1, 1, 1, 0, 0, 0, 0, 0, 0, 0, 
0, 0 )_{({\sf diagonal})}.
\label{eq:PdgSO4}
\end{eqnarray}
The prescription for general quantized orthogonal group 
$SO_q ( N )$ is evident.

A transformation that diagonalizes the generalized braid matrix
${\hat R} ( v, w )$ must diagonalize each projector {\em separately},
as these projectors are functions of the matrix ${\hat R} ( v, w )$.
Using the above form of the diagonalized projectors, we avoid 
introducing the inverse of the diagonalizing matrix $M$ as follows.
Assuming that a diagonalizing matrix of non-zero determinant exist,
we have
\begin{equation}
M P_{(0)} = P_{(0)}^{({\sf d})} M,\qquad M P_{(-)} = 
P_{(-)}^{({\sf d})} M,	
\label{eq:MPPd}
\end{equation}
where the diagonalized projectors $P^{({\sf d})}$ have explicit 
forms given by (\ref{eq:PdgSO3}) and (\ref{eq:PdgSO4}), and their 
evident generalizations. The projectors $P_{(0)}$ and $P_{(-)}$ 
being known {}from the standard results, we avoid direct 
introduction of the matrix $M^{-1}$, which is nonlinear in 
elements of $M$. 

The diagonalizing relations (\ref{eq:MPPd}) generate  a set of 
{\em linear} constraints on the elements of $M$. The coefficients 
in each equation are fully known. They will in general leave room 
for many arbitrary parameters in the matrix $M$, subject to the 
constraint that it has a non-zero determinant. This arbitrariness 
can be `factored out' as follows. Supposing that in the case of the 
the quantum group $SO_q (3)$ we have found a convenient solution 
${\hat M}$ for the diagonalizing matrix, the relation 
(\ref{eq:MPPd}) in conjunction with (\ref{eq:PdgSO3}) yield 
the structure
\begin{eqnarray}
{\hat M} {\hat R} ( v, w ) {\hat M}^{-1} &=& {\hat M} \left( 
{\sf I} + v P_{(-)} + w P_{(0)} \right) {\hat M}^{-1}\nonumber\\   
\phantom{{\hat M} {\hat R} ( v, w ) {hat M}^{-1}} &=& ( 1 + w, 
1 + v, 1 + v, 1 + v, 1, 1, 1, 1, 1 )_{( {\sf diagonal})}.
\label{eq:Mhat}
\end{eqnarray} 
Keeping the above block structure in mind we now use a block 
diagonal matrix 
\begin{equation}
{\sf M} = \left( {\sf M}_{(0)}, {\sf M}_{(-)}, {\sf M}_{(+)}
\right)_{( {\sf block\,diagonal} )},
\label{eq:UMblock}
\end{equation}  
where ${\sf M}_{(0)}, {\sf M}_{(-)}$ and ${\sf M}_{(+)}$ 
respectively are  $1 \times 1, 3 \times 3$ and $5 \times 5$
invertible matrices with non-zero determinant. Using the 
property (\ref{eq:Mhat}) and the above structure of the matrix 
${\sf M}$, we  obtain
\begin{equation}
\left( {\sf M} {\hat M} \right) {\hat R} ( v, w )  
\left( {\sf M} {\hat M} \right)^{-1} = {\hat M} {\hat R} ( v, w ) 
{\hat M}^{-1}.
\label{eq:MMR} 
\end{equation}
Thus the diagonalization is preserved exactly. The above matrix 
${\sf M}_{(0)}$ is a non-zero constant, whereas the matrices
${\sf M}_{( \mp )}$, apart from the single constraint of 
invertibility have {\em arbitrary} elements. This is the source 
of arbitrariness in the diagonalization procedure. For the general 
case the dimensions of ${\sf M}_{(0)}, {\sf M}_{(-)}$ and 
${\sf M}_{(+)}$ are given by the respective traces 
(\ref{eq:Trexpl}) of the corresponding projectors.

But how do we select ${\hat M}$? For the quantum groups $SO_q (3)$
and $SO_q (4)$ we present below the diagonalizing matrix ${\hat M}$
possessing a remarkably helpful feature: {\em the rows of ${\hat M}$
are mutually orthogonal}. It is presumably possible to choose the 
diagonalizing matrix ${\hat M}$ retaining this feature for all 
orthogonal quantum groups. But here we do not attempt to construct
such a general solution for an arbitrary case. 

For the quantum group $SO_q (3)$ we define the parameters
$s = - q^{- \half} ( 1 - q ), t = - q^{- \frac{3}{2}} ( 1 + q )$,
and obtain the relevant diagonalizing matrix:
\begin{equation}
{\hat M } = \left(
\begin{array}{ccccccccc}
0 & 0 & 1 & 0 & q^{\half} & 0 & q & 0 & 0\\
0 & 1 & 0 & - q & 0 & 0 & 0 & 0 & 0\\
0 & 0 & 0 & 0 & 0 & 1 & 0 & - q & 0\\
0 & 0 & 1 & 0 & s & 0 & -1 & 0 & 0\\
1 & 0 & 0 & 0 & 0 & 0 & 0 & 0 & 0\\
0 & 1 & 0 & q^{-1} & 0 & 0 & 0 & 0 & 0\\ 
0 & 0 & 1 & 0 & t & 0 & q^{-2} & 0 & 0\\
0 & 0 & 0 & 0 & 0 & 1 & 0 & q^{-1} & 0\\
0 & 0 & 0 & 0 & 0 & 0 & 0 & 0 & 1\\
\end{array}
\right).
\label{eq:MhatSO3}
\end{equation}
Here we note the existence of the following orthogonal triplets:
\begin{equation}
( 1, q^{\half}, q ), ( 1, s, - 1 ), ( 1, t, q^{-2} ),
\label{eq:triplet}
\end{equation}
which, as will be seen later, play a role in the context of the 
noncommutative coordinates. For the quantum group $SO_q (4)$, 
one possible choice of the ordering of the rows gives the following 
diagonalizing matrix:  
\begin{equation}
{\hat M} = \left(
\begin{array}{cccccccccccccccc}
0 & 0 & 0 & 1 & 0 & 0 & q & 0 & 0 & q & 0 & 0 & q^2 & 0 & 0 & 0\\
0 & 0 & 0 & 1 & 0 & 0 & q & 0 & 0 & - q^{- 1} & 0 & 0 & -1 & 0 &
0 & 0\\
0 & 1 & 0 & 0 & - q & 0 & 0 & 0 & 0 & 0 & 0 & 0 & 0 & 0 & 0 & 0\\
0 & 0 & 1 & 0 & 0 & 0 & 0 & 0 & - q & 0 & 0 & 0 & 0 & 0 & 0 & 0\\
0 & 0 & 0 & 0 & 0 & 0 & 0 & 1 & 0 & 0 & 0 & 0 & 0 & - q & 0 & 0\\
0 & 0 & 0 & 0 & 0 & 0 & 0 & 0 & 0 & 0 & 0 & 1 & 0 & 0 & - q & 0\\
0 & 0 & 0 & 1 & 0 & 0 & - q^{- 1} & 0 & 0 & q & 0 & 0 & - 1 & 0 &
0 & 0\\
1 & 0 & 0 & 0 & 0 & 0 & 0 & 0 & 0 & 0 & 0 & 0 & 0 & 0 & 0 & 0\\
0 & 1 & 0 & 0 & q^{-1} & 0 & 0 & 0 & 0 & 0 & 0 & 0 & 0 & 0 & 0 &
0\\
0 & 0 & 1 & 0 & 0 & 0 & 0 & 0 & q^{-1} & 0 & 0 & 0 & 0 & 0 & 0 &
0\\
0 & 0 & 0 & 1 & 0 & 0 & - q^{- 1} & 0 & 0 & - q^{- 1} & 0 & 0 &
q^{- 2} & 0 & 0 & 0\\
0 & 0 & 0 & 0 & 0 & 0 & 0 & 1 & 0 & 0 & 0 & 0 & 0 & q^{- 1} & 0 &
0\\
0 & 0 & 0 & 0 & 0 & 0 & 0 & 0 & 0 & 0 & 0 & 1 & 0 & 0 & q^{-1} & 
0\\
0 & 0 & 0 & 0 & 0 & 0 & 0 & 0 & 0 & 0 & 1 & 0 & 0 & 0 & 0 & 0\\
0 & 0 & 0 & 0 & 0 & 1 & 0 & 0 & 0 & 0 & 0 & 0 & 0 & 0 & 0 & 0\\
0 & 0 & 0 & 0 & 0 & 0 & 0 & 0 & 0 & 0 & 0 & 0 & 0 & 0 & 0 & 1\\
\end{array}
\right).
\label{eq:MhatSO4}
\end{equation}
Here the relevant orthogonal quadruplets read
\begin{equation}
( 1, q, q, q^2 ), ( 1, q, - q^{-1}, - 1 ), ( 1, - q^{-1}, q, - 1 ),
( 1, - q^{-1}, - q^{- 1}, q^{-2} ).
\label{eq:quadruplet}
\end{equation}   

It may be noted that the orthogonality of the rows delivers the inverse
of the diagonalizing matrix ${\hat M}^{-1}$ effortlessly. To this end, 
we take the transposed matrix ${\hat M}^{\sf T}$ and normalize each
element of its column $j$ by the same factor $c_j$ so that the sum
of squares of all elements of this column $j$ multiplied with $c_j$
equates to unity. For the diagonalizing matrix ${\hat M}$ in
(\ref{eq:MhatSO3}), these normalization constants are given by
\begin{eqnarray}    
&&q ( 1 + [2] )\,c_1 = q [2]\,c_2 = q [2]\,c_3 = [2]\,c_4 = c_5
= q^{-1} [2]\,c_6\nonumber\\ 
&&= q^{-2} [2] ( 1 + [2] )\,c_7 = q^{-1} [2]\,c_8 = c_9 = 1,
\label{eq:SO3c}
\end{eqnarray}  
whereas the corresponding normalization constants for the matrix
in (\ref{eq:MhatSO4}) read
\begin{eqnarray}
&&q^2 [2]^2\,c_1 = [2]^2\,c_2 = q [2]\,c_3 = q [2]\,c_4 
= q [2]\,c_5 = q [2]\,c_6 = [2]^2\,c_7 = c_8 = q^{-1} [2]\,c_9 
\nonumber\\
&&= q^{-1} [2]\,c_{10} = q^{-2} [2]^2\,c_{11} = q^{-1} [2]\,c_{12} 
= q^{-1} [2]\,c_{13} = c_{14} = c_{15} = c_{16} = 1.
\label{eq:SO4c}
\end{eqnarray}

Now we show that the orthogonality of the rows also directly leads 
to the eigenvectors of the generalized braid matrix 
${\hat R} ( v, w )$. We again choose the quantum group $SO_q (3)$ 
as a typical example. Let 
$\{ |V_{k}\rangle\,| ( k = 1,\cdots,9 ) \}$ be the eigenvectors
of the ${\hat R} ( v, w )$ matrix with the eigenvalues $a_k$:
\begin{equation}
{\hat R} ( v, w )\,|V_k\rangle = a_k\,|V_k\rangle,
\label{eq:eigenvec} 
\end{equation}
where no summation over the index $k$ on the right is implied. 
Multiplying the equation (\ref{eq:eigenvec}) with the matrix 
${\hat M}$ on both sides and employing diagonalization property 
(\ref{eq:Mhat}) we obtain
\begin{equation}
( 1 + w, 1 + v, 1 + v, 1 + v, 1, 1, 1, 1, 1 )_{( {\sf diagonal} )} 
| {\hat M}\,V_k\rangle = a_k |{\hat M}\,V_k\rangle.
\label{eq:MVk}
\end{equation} 
Hence, apart {}from overall constant factors, we may choose the 
transposed vectors as follows:
\begin{equation}
{| {\hat M} \,V_k \rangle}^{\sf T} = \{ ( 1, 0,\cdots, 0 ),\,
( 0, 1, 0,\cdots, 0 ),\,\cdots,( 0,\cdots,0, 1 ) \}.
\label{eq:MVvec}
\end{equation}
The eigenvalues are given by the diagonal elements of the matrix
${\hat M} {\hat R} ( v, w ) {\hat M}^{-1}$ as exhibited in
(\ref{eq:Mhat}). Hence the obvious construction $| V_k \rangle
\, = {\hat M}^{- 1}\,| {\hat M}\,V_k \rangle$,  
along with (\ref{eq:MVvec}) furnishes the eigenvector 
$| V_k \rangle$ as the $k$-th column of ${\hat M}^{-1}$. Now, in
view of the construction of the inverse matrix ${\hat M}^{-1}$ 
previously discussed preceding (\ref{eq:SO3c}), the eigenvector
$| V_k \rangle$ {\em is finally given by the $k$-th row of the
matrix} ${\hat M}$. Replacing here the matrix ${\hat M }$ by its
alternatives ${\sf M} {\hat M}$ amounts to, as given in 
(\ref{eq:MMR}), taking linear combinations of the above 
eigenvectors with the same  eigenvalue. The eigenvectors and 
particularly the highest eigenvalue, dependent here on the choice 
of the spectral variables $( v, w )$, of the generalized braid 
matrix ${\hat R} ( v, w )$ are of interest in related models
of statistical mechanics. Another, quite different, interest 
in the above diagonalization of the ${\hat R} ( v, w )$ matrix
will be pointed out in Sec. 6 in the context of noncommutative
spaces.  

So far we have considered the orthogonal quantum group 
$SO_q ( N )$. A parallel, but much simpler, formalism may be 
developed for the linear quantum group $GL_q ( N )$, as only two 
projectors are present there. Instead of discussing the general 
case, we, for the purpose of illustration, consider a 
biparametric $( p, q )$ deformation of the group  $GL (2)$ in the 
remaining part of the present section. In \cite{C01} one of us 
introduced the generalized braid matrix ( for $p \neq 0$ )
\begin{equation}
{\hat R} ( K; p, q ) =
\left(
\begin{array}{cccc}
1 & 0 & 0 & 0\\
0 & ( 1 - K ) & K p^{- 1} & 0\\
0 & K q & ( 1 - K q p^{-1} ) & 0\\
0 & 0 & 0 & 1\\
\end{array}
\right),
\label{eq:Rpq}
\end{equation}
which satisfies the strict braid equation for the two parametric
values $K = 1, p q^{-1}$. Maintaining a specific parametrization 
of the spectral variable $v$, we write in conformity with the 
notations of the Sec. 2:
\begin{equation}
{\hat R} ( K; p, q ) = {\sf I} - K ( 1 + q p^{-1} )\,P_{(-)},
\label{eq:Rpqproj}
\end{equation}
where the projector is
\begin{equation}
P_{(-)} = {( 1 + q p^{-1} )}^{-1}
\left(                                                    
\begin{array}{cccc}
0 & 0 & 0 & 0\\
0 & 1 & - p^{-1} & 0\\
0 & - q & q p^{-1} & 0\\
0 & 0 & 0 & 0\\                 
\end{array}
\right).
\label{eq:Ppq}
\end{equation}
The projector $P_{(-)}$ may be diagonalized by
\begin{equation}
{\hat M} = 
\left(
\begin{array}{cccc}
0 & - 1 & p^{-1} & 0\\
0 & q & 1 & 0\\
1 & 0 & 0 & 1\\
1 & 0 & 0 & - 1\\
\end{array}
\right)
\label{eq:Mpq} 
\end{equation}
giving ${\hat M} P_{(-)} {\hat M}^{-1} = {( 1, 0, 0, 
0 )}_{( {\sf diagonal})}$. This, following the spectral 
decomposition (\ref{eq:Rpqproj}), provides the diagonalized 
${\hat R}$ matrix. In this simple case the matrix 
${\hat M}^{-1}$ may be easily computed for all values of 
deformation parameters $( p, q )$. But except for the special 
case $p q = 1$, the first two rows of the present diagonalizer 
${\hat M}$ are not orthogonal to each other.

\sect{Generalized braid matrices and consequent noncommutative 
spaces}

We first present the general prescription for covariant 
quantization of spaces implementing the generalized braid matrices
${\hat R} ( v )$ and ${\hat R} ( v, w )$. Our prescription follow
the standard structure except for the presence of the arbitrary
values of the spectral variables $( v, w )$. Then we will display 
how the diagonalizations of the generalized braid matrices
presented in Sec. 5 enable us to extract in a convenient fashion
the contents of the above prescription using, as before, the 
quantum groups $GL_q ( 2 ), SO_q ( 3 )$ and $SO_q ( 4 )$ as 
examples. Generalization to higher dimensional quantum groups can 
then be easily carried out. Our study will thus be limited. 
We hope to present elsewhere a fuller exploration of the possible 
roles of the variables $( v, w )$. 
 
For the quantum group $GL_q (2)$, we use the projector 
(\ref{eq:Ppq}) with the restriction $p = q^{-1}$, and set the
${\hat R} ( v )$ matrix as in (\ref{eq:ARv}) satisfying the
characteristic equation (\ref{eq:Afact}). Using the standard 
notations with coordinates and differentials given, respectively,
by $\{ x_i, \xi_i \equiv dx_i\,|\, i = ( 1, 2 ) \}$, we obtain 
\begin{equation}
P_{(-)} x \otimes x = 0
\label{eq:PMxx}
\end{equation}
and
\begin{equation}
x \otimes \xi = B\,\xi \otimes x,
\label{eq:xxi}
\end{equation}
whereas Leibnitz rule and covariance lead to
\begin{equation}
( B + {\sf I} )\,\xi \otimes \xi = 0.
\label{eq:xixi}
\end{equation}
The matrix $B$ reads
\begin{equation}
B = - {\sf I} + \mu\,\left( {\hat R} ( v ) - ( 1 + v ) 
{\sf I}\right), 
\label{eq:BR}
\end{equation}
while the parameter $\mu$ is arbitrary. We also note that the 
above structure of $B$ ensures the orthogonality: $( B + {\sf I} ) 
P_{(-)} = 0$. Using the diagonalizer ${\hat M}$ given in 
(\ref{eq:Mpq}), and
adapting it to the parametric choice $p = q^{-1}$, we evaluate
\begin{equation}
{\hat M} x \otimes \xi = {\hat M} 
\left(  
\begin{array}{c}
x_1 \xi_1\\
x_1 \xi_2\\
x_2 \xi_1\\
x_2 \xi_2\\
\end{array}
\right) = \left(
\begin{array}{c}
- x_1 \xi_2 + q x_2 \xi_1\\
q x_1 \xi_2 + x_2 \xi_1\\
x_1 \xi_1 +  x_2 \xi_2\\
x_1 \xi_1 -  x_2 \xi_2\\
\end{array}
\right).
\label{eq:Mxxi}
\end{equation}
The analogous results for ${\hat M} x \otimes x,\,{\hat M} \xi 
\otimes \xi$ and ${\hat M} \xi \otimes x$ can be read off 
(\ref{eq:Mxxi}) readily. Using the previous explicit constructions
we may now obtain the following diagonalized structures
\begin{eqnarray}
&&{\hat M} P_{(-)} {\hat M}^{-1} = {( 1, 0, 0, 
0 )}_{({\sf diagonal})},\nonumber\\
&&{\hat M} {\hat R} ( v ) {\hat M}^{-1} = 
{( 1 + v, 1, 1, 1 )}_{({\sf diagonal})},\nonumber\\
&&{\hat M} ( B + {\sf I} ) {\hat M}^{-1} = 
\mu {( 0, - v, - v, - v )}_{({\sf diagonal})},\nonumber\\
&&{\hat M} B {\hat M}^{-1} = 
- {( 1, 1 + \mu v, 1 + \mu v, 1 + \mu v )}_{({\sf diagonal})},
\label{eq:GL2DG}
\end{eqnarray}
which in conjunction with (\ref{eq:Mxxi}) allow us to immediately
write down the relations
\begin{equation}
x_1 x_2 = q\, x_2 x_1,\quad \xi_1 \xi_2 = - q^{- 1}\xi_2 \xi_1,
\quad\xi_1^2 = 0,\quad \xi_2^2 = 0
\label{eq:Manin}
\end{equation}
and
\begin{eqnarray}
&&( - x_1 \xi_2 + q x_2 \xi_1 ) =  
- ( - \xi_1 x_2 + q \xi_2 x_1 ),\nonumber\\
&&( q x_1 \xi_2 + x_2 \xi_1 ) = 
- ( 1 + \mu v ) ( q \xi_1 x_2 + \xi_2 x_1 ),\nonumber\\
&&( x_1 \xi_1 + x_2 \xi_2 ) = 
- ( 1 + \mu v ) ( \xi_1 x_1 + \xi_2 x_2 ),\nonumber\\
&&( x_1 \xi_1 -  x_2 \xi_2 ) = 
- ( 1 + \mu v ) ( \xi_1 x_1 - \xi_2 x_2 ).
\label{eq:xxibrd}
\end{eqnarray}
Adapting notations and parametrizations the corresponding results 
of \cite{C01} may be obtained. But the consequences of the 
diagonalization is particularly striking in the modular structure 
evident in (\ref{eq:xxibrd}). Linear combinations of $x_i \xi_j$ 
are picked out that are {\em proportional to the same combination} 
of $\xi_i x_j$ on the right. This is, as will be seen, a general 
feature in all cases. It was shown in \cite{C01} that 
\begin{equation}
\Phi^2 =  ( \xi_1 x_2 - q \xi_2 x_1 )^2 = 0.
\label{eq:Phinil}
\end{equation} 
Here, using (\ref{eq:xxibrd}) it immediately follows that the 
combination $( x_1 \xi_2 - q x_2 \xi_1 )$ is  also nilpotent. Our 
formalism also signals clearly special values of the parameters. 
It is evident {}from (\ref{eq:xxibrd}) that $\mu = - v^{-1}$ is 
a very special case. For the quantum group $GL_q ( N )$ for an 
arbitrary $N$ the prescriptions (\ref{eq:PMxx}) to (\ref{eq:BR}) 
remain the same except that higher dimensional $N^2 \times N^2$ 
matrices need to be considered.

Now we consider the noncommuting spaces associated with the
generalized braid matrix ${\hat R} ( v, w )$ of the orthogonal 
quantum group $SO_q ( N )$. This matrix has the structure 
(\ref{eq:BCDRvw}) and it satisfies the characteristic equation 
(\ref{eq:BCDfact}). The braiding structure (\ref{eq:PMxx}) of the 
noncommutative coordinates   now may also be expressed as
\begin{equation}
\left( {\hat R} ( v, w ) - {\sf I} \right)
\left( {\hat R} ( v, w ) - ( 1 + w )\,{\sf I} \right)
x \otimes x = 0.
\label{eq:RRxx} 
\end{equation}
Leibnitz rule and covariance ensure equations (\ref{eq:xxi})
and (\ref{eq:xixi}) with the matrix $B$ given by (\ref{eq:BR}).
As there are three projectors for the orthogonal quantum groups,
the braiding relation (\ref{eq:xixi}) now reduces to
\begin{equation}
P_{(+)} \xi \otimes \xi = 0,\quad P_{(0)} \xi \otimes \xi = 0. 
\label{eq:Pxixi}
\end{equation}

Now we will explicitly demonstrate the above structure for the 
quantum groups $SO_q (3)$ and $SO_q (4)$. For the quantum group
$SO_q (3)$ the diagonalizing matrix ${\hat M}$ is given in  
(\ref{eq:MhatSO3}). The diagonalized operators now read
\begin{eqnarray}
&&{\hat M} P_{(-)} {\hat M}^{-1} = {( 0, 1, 1, 1, 0, 0, 0, 
0, 0  )}_{({\sf diagonal})},\nonumber\\
&&{\hat M} {\hat R} ( v, w ) {\hat M}^{-1} = {( 1 + w, 1 + v, 
1 + v, 1 + v, 1, 1, 1, 1, 1 )}_{({\sf diagonal})},\nonumber\\
&&{\hat M} B {\hat M}^{-1} = - ( 1 + \mu ( v- w ), 1, 1, 1,
1 + \mu v, 1 + \mu v, 1 + \mu v, 
1 + \mu v, 1 + \mu v )_{({\sf diagonal})}.
\label{eq:SO3DG}
\end{eqnarray}
Here we choose, as in \cite{C02}, the triplets $( x_{-}, x_{0}, 
x_{+} )$ and $( \xi_{-}, \xi_{0}, \xi_{+} )$ as the basis elements
for the noncommuting coordinates and the differentials 
respectively. Using the diagonalizing matrix ${\hat M}$ given
in (\ref{eq:MhatSO3}) we now compute
\begin{equation}
{\hat M} x \otimes \xi = {\hat M} \left(
\begin{array}{c}
x_{-} \xi_{-}\\
x_{-} \xi_{0}\\
x_{-} \xi_{+}\\
x_{0} \xi_{-}\\
x_{0} \xi_{0}\\
x_{0} \xi_{+}\\
x_{+} \xi_{-}\\
x_{+} \xi_{0}\\
x_{+} \xi_{+}\\
\end{array}
\right) = \left(
\begin{array}{c}
x_{-} \xi_{+} + q^{1/2} x_{0} \xi_{0} + q x_{+} \xi_{-}\\
x_{-} \xi_{0} - q x_{0} \xi_{-}\\
x_{0} \xi_{+} - q x_{+} \xi_{0}\\
x_{-} \xi_{+} + s x_{0} \xi_{0} - x_{+} \xi_{-}\\
x_{-} \xi_{-}\\
x_{-} \xi_{0} + q^{-1} x_{0} \xi_{-}\\
x_{-} \xi_{+} + t x_{0} \xi_{0} + q^{-2} x_{+} \xi_{-}\\
x_{0} \xi_{+} + q^{-1} x_{+} \xi_{0}\\
x_{+} \xi_{+}\\
\end{array}
\right).
\label{eq:MxxiSO3}
\end{equation}
The matrices  ${\hat M} x \otimes x, {\hat M} \xi \otimes \xi$ and
${\hat M} \xi \otimes x$ have evident analogous forms. The 
parameters $( \mu, v, w )$ do not appear in the braiding
structures (\ref{eq:PMxx}) and (\ref{eq:Pxixi}). Therefore these 
braiding relations have the usual form, given, for example, in 
Eq. (3.48) of \cite{C02}. For the sake of completeness we present 
them here:
\begin{eqnarray}
&&x_{-} x_0 = q x_0 x_{-},\quad x_0 x_{+} = q x_{+} x_0,\quad 
x_{+} x_{-} - x_{-} x_{+} = ( q^{1/2} - q^{- 1/2} ) x_{0}^2,
\nonumber\\ 
&&\xi_{-}^2 = 0,\qquad \xi_{+}^2 = 0,\qquad
\xi_{-} \xi_{+} + \xi_{+} \xi_{-} = 0,\nonumber\\
&&q \xi_{-} \xi_{0} + \xi_{0} \xi_{-} = 0,\quad
q \xi_{0} \xi_{+} + \xi_{+} \xi_{0} = 0,\quad
\xi_{0}^2 = ( q^{1/2} - q^{- 1/2} ) \xi_{-} \xi_{+},
\label{eq:stdbrSO3} 
\end{eqnarray} 
The constraints due to (\ref{eq:xxi}) do involve the parameters
$( \mu, v, w )$ and have been obtained in Eq. (3.48) of \cite{C02}. 
Here we present them in the form directly given by (\ref{eq:SO3DG})
and (\ref{eq:MxxiSO3}). After implementing our diagonalization,
we obtain
\begin{eqnarray}
( x_{-} \xi_{+} + q^{1/2} x_{0} \xi_{0} + q x_{+} \xi_{-} ) &=&
- (1 + \mu ( v - w ))
( \xi_{-} x_{+} + q^{1/2} \xi_{0} x_{0} + q\,\xi_{+} x_{-} ),
\nonumber\\
( x_{-} \xi_{+} + s x_{0} \xi_{0} - x_{+} \xi_{-} ) &=&
- ( \xi_{-} x_{+} + s\,\xi_{0} x_{0} - \xi_{+} x_{-} ),\nonumber\\
( x_{-} \xi_{+} + t x_{0} \xi_{0} + q^{-2} x_{+} \xi_{-} ) &=&
- ( 1 + \mu v ) 
( \xi_{-} x_{+} + t\,\xi_{0} x_{0} + q^{-2} \xi_{+} x_{-} ),
\nonumber\\
( x_{-} \xi_{0} - q x_{0} \xi_{-} ) &=&
- ( \xi_{-} x_{0} - q\,\xi_{0} x_{-} ),\nonumber\\
( x_{-} \xi_{0} + q^{-1} x_{0} \xi_{-} ) &=&
- ( 1 + \mu v ) ( \xi_{-} x_{0} + q^{-1} \xi_{0} x_{-} ),
\nonumber\\
( x_{0} \xi_{+} - q x_{+} \xi_{0} ) &=&
- ( \xi_{0} x_{+} - q \xi_{+} x_{0} ),\nonumber\\
( x_{0} \xi_{+} + q^{-1} x_{+} \xi_{0} ) &=&
- ( 1 + \mu v ) ( \xi_{0} x_{+} + q^{-1} \xi_{+} x_{0} ),
\nonumber\\
x_{-} \xi_{-} &=& - ( 1 + \mu v ) \xi_{-} x_{-},\nonumber\\
x_{+} \xi_{+} &=& - ( 1 + \mu v ) \xi_{+} x_{+}.
\label{eq:modbrSO3}
\end{eqnarray}
Thus, as signalled before, the diagonalization selects out linear 
combinations, which are proportional under the operation
$x_i \xi_j \rightarrow \xi_i x_j$. We note that the coefficients 
appearing in the triplets, namely: $( 1, q^{1/2}, q ),
( 1, - q^{- 1/2} ( 1 - q ), - 1 ), ( 1, - q^{- 3/2} ( 1 + q ),
q^{-2} )$ are mutually orthogonal. The same property is evident
with the doublets: $( 1, - q ), ( 1, q^{-1} )$. The parametric
values $\mu = - v^{-1}$ again generate a special case. Using 
(\ref{eq:modbrSO3}) each  $x_i \xi_j$ may be written as linear
combination of $\xi_k x_l$ terms. This is provided in \cite{C02}.
Here we want to emphasize the modular structure of 
(\ref{eq:modbrSO3}). In the preceding discussion regarding the 
quantum group $GL_q (2)$ it was noted that a nilpotent bilinear 
structure, given in (\ref{eq:Phinil}), arise directly out of the 
diagonalization process. Presence of such structures in 
(\ref{eq:modbrSO3}) should be sought. A detailed study of our 
generalized spaces, arising as a consequence of generalized braid
matrices, will be presented elsewhere.

To investigate our noncommutative spaces associated with the
quantum group $SO_q (4)$, we proceed exactly as in the previous 
example. Now the diagonalizing matrix ${\hat M}$ is given by
(\ref{eq:MhatSO4}). The basis for the coordinates and the 
corresponding basis for the differentials are denoted by
$( x_1, x_2, x_3, x_4 )$ and $( \xi_1, \xi_2, \xi_3, \xi_4 )$ 
respectively. The braiding relations for the bilinears $x \otimes 
x$ and $\xi \otimes \xi$ are independent of the parameters 
$( \mu, v, w )$ introduced here; and, consequently, remain  
the standard noncommutativity constraints. We exhibit below,
for brevity, only the relevant modular structure analogous to
(\ref{eq:modbrSO3}). Now we have the quadruplets, doublets and 
singlets with typically orthogonalized constraints as before:
\begin{eqnarray}
( x_1 \xi_4 + q x_2 \xi_3 + q x_3 \xi_2 + q^2 x_4 \xi_1 ) &=&
- (1 + \mu ( v - w ))
( \xi_1 x_4 + q \xi_2 x_3 + q \xi_3 x_2 + q^2 \xi_4 x_1 ),
\nonumber\\
( x_1 \xi_4 + q x_2 \xi_3 - q^{-1} x_3 \xi_2 - x_4 \xi_1 ) &=&
- ( \xi_1 x_4 + q \xi_2 x_3 - q^{-1} \xi_3 x_2 - \xi_4 x_1 ),
\nonumber\\
( x_1 \xi_4 - q^{-1} x_2 \xi_3 + q x_3 \xi_2 - x_4 \xi_1 ) &=&
- ( \xi_1 x_4 - q^{-1} \xi_2 x_3 + q \xi_3 x_2 - \xi_4 x_1 ),
\nonumber\\    
( x_1 \xi_4 - q^{-1} x_2 \xi_3 - q^{-1} x_3 \xi_2
+q^{-2} x_4 \xi_1 ) &=& - ( 1 + \mu v ) ( \xi_1 x_4 
- q^{-1} \xi_2 x_3 - q^{-1} \xi_3 x_2 + q^{-2} \xi_4 x_1 ),
\nonumber\\
( x_i \xi_j - q x_j \xi_i ) &=&
- ( \xi_i x_j - q \xi_j x_i ),\nonumber\\
( x_i \xi_j + q^{-1} x_j \xi_i ) &=&
- ( 1 + \mu v ) ( \xi_i x_j + q^{-1} \xi_j x_i ),\nonumber\\
x_k \xi_k &=& - ( 1 + \mu v ) \xi_k x_k,
\label{eq:modbrSO4}
\end{eqnarray}
where $( i, j ) = \{ ( 1, 2 ), ( 1, 3 ), ( 2, 4 ), 
( 3, 4 ) \}$ and $k = \{ 1, 2, 3, 4 \}$. No summation over the 
index $k$ is meant in the last equation in (\ref{eq:modbrSO4}).
Comments parallel to those following (\ref{eq:modbrSO3}) are 
relevant here.

Here we sum up the above procedure for the quantum group 
$SO_q (N)$ with an arbitrary value of $N$. The formalism 
described in equations (\ref{eq:PMxx})-(\ref{eq:BR}), and also in 
(\ref{eq:Pxixi}) holds for an arbitrary value of $N$. We have
emphasized on the modular structure and the parameter dependence
of the braiding relations arising {}from (\ref{eq:xxi}) and
our diagonalization of the generalized braid matrices. Setting the
spectral variables $( v, w )$ equal to their braid values
(\ref{eq:brdval}), and also fixing the parameter $\mu = q^2$, we
may recover the standard quantization prescriptions. So far we 
have implicitly assumed that the spectral variables are 
nonvanishing: $v \neq 0, w \neq 0$. When one of them vanishes, 
the characteristic equation obeyed by the generalized braid 
matrix becomes, as discussed for instance in (\ref{eq:Rwchar}), 
quadratic rather than cubic. The corresponding reformulation of 
the prescription for our noncommutative spaces is straightforward 
and will not be presented here.

\sect{Remarks}
The following braid matrices,
\begin{equation}
{\hat R} = \left(
\begin{array}{cccc}
1 & 0 & 0 & 1\\
0 & 1 & - 1 & 0\\
0 & 1 & 1 & 0\\
- 1 & 0 & 0 & 1\\
\end{array}
\right)
\label{eq:S03}
\end{equation}
and
\begin{equation}
{\hat R} = \left(
\begin{array}{cccc}
0 & 0 & 0 & q\\
0 & 1 & 0 & 0\\
0 & 0 & 1 & 0\\
q & 0 & 0 & 0\\
\end{array}
\right),
\label{eq:S14}
\end{equation}
were studied in detail \cite{ACDM02} in the context of exotic 
bialgebras \cite{ACDM01,ACDM02}. These authors investigated  the 
spectral decomposition of these braid matrices and the 
Baxterization thereof. These respective cases were denominated 
as $S03$ and $S14$ in the classification scheme of $4 \times 4$ 
braid matrices presented in \cite{H93}. The analogies and 
differences of the ${\hat R}$ matrix in (\ref{eq:S03}) with the 
unitary case, and that in (\ref{eq:S14}) with the orthogonal case 
have been discussed \cite{ACDM01, ACDM02} earlier. Noncommutative 
spaces associated with these two and other `exotic' $4 \times 4$ 
${\hat R}$ matrices have been studied \cite{ACDM01, ACDM02}. These 
$4 \times 4$ matrices classified in \cite{H93} are `exotic' in 
the sense that they are not obtained by restricting some universal 
${\hat{\cal R}}$ matrix to this dimension. They represent distinct 
supplementary possibilities for the dimension $4 \times 4$.

How do we construct higher dimensional analogues of these matrices? 
In Sec. 4 we have presented a canonical construction of a class of 
$N^2 \times N^2$ braid matrices complete with the corresponding 
MBE and the Baxterized forms of these matrices. In fact for each 
even $N$ we have {\em two} solutions: `exotic orthogonal' for the 
choice $\varepsilon = 1$, and `exotic symplectic' for  the 
parametric value $\varepsilon = - 1$. But these class presumably 
does not exhaust such possibilities for each $N$. Our general 
reduction of trilinear terms in Sec. 3 insistently pointed out 
the class studied in Sec. 4. How do we investigate other 
possibilities? Generalization of Hietarinta's approach \cite{H93} 
to higher dimensional cases would be extremely laborious. Still 
a more thorough search may be worthwhile.

Several interesting aspects of our formalism have not been 
addressed in the present work. Quasi-classical limits, 
$L$ operators and Yangians are relevant examples. Applications,
particularly of our new class of solutions, to integrable models
\cite{R85, B85} would be worth exploring. Certain specific 
properties would be lost if the parameters $( \mu, v, w )$ move
away from their standard values. We wish to study
new interesting features which may emerge for other values of      
the parameters $( \mu, v, w )$. The present approach via 
diagonalization may be helpful.
      
Let us end by taking a closer look at the mutually orthogonal 
sets of triplets and quadruplets, appearing in (\ref{eq:triplet}) 
and (\ref{eq:quadruplet}) as a consequence of our 
diagonalizations. For the quantum group $SO_q (3)$ the constraint 
(\ref{eq:PMxx}) contain, corresponding to the set 
$( 1, ( q^{\half} - q^{- \half} ), - 1 )$ in (\ref{eq:triplet}), 
the commutation relation, listed before in the set 
(\ref{eq:stdbrSO3}) but remodeled here for the purpose 
of convenience: 
\begin{equation} 
x_{-} x_{+} - x_{+} x_{-} + ( q^{\half} -
q^{- \half} ) {x_{0}}^2 = 0. 
\label{eq:xcomm} 
\end{equation} 
A constraint trivially true in the commutative limit $( q = 1 )$ 
is thus consistently maintained. The analogous expressions 
corresponding to the other two triplets in the set 
(\ref{eq:triplet}), selected out by the other two diagonalized 
projectors $P_{(0)}$ and $P_{(+)}$ are {\em not} constrained
to be zero. Introducing the metric and the star operation 
( Ex. 4.1.22 in \cite{M99} ) as 
\begin{equation} 
\left( x_{\pm} \right)^{*} = q^{\mp \half} \left( x_{\mp} \right), 
\qquad {x_{0}}^{*} = x_{0},
\label{eq:xstar} 
\end{equation} 
the other two triplets mentioned above lead to the surfaces with 
invariants ${\sf k}_1$ and ${\sf k}_2$:
\begin{eqnarray} 
&&{x_{-}}^{*} x_{-} + {x_{+}}^{*} x_{+} + {x_{0}}^{*}
x_{0} = {\sf k}_1,\nonumber\\ 
&&q^{- \frac{3}{2}} {x_{-}}^{*} x_{-} + q^{\frac{3}{2}} 
{x_{+}}^{*} x_{+} - ( q^{\half} + q^{- \half} ) 
\,{x_{0}}^{*} x_{0} = {\sf k}_2, 
\label{eq:NCsurface}
\end{eqnarray}
where ${\sf k}_1$ is usually denoted as the `distance squared'
{\em i.e.} ${\sf k}_1 \equiv {\sf r}^2 \geq 0$. 
The above two surfaces denote a $q$-deformed sphere and a 
$q$-deformed hyperboloid respectively. In the context of our 
diagonalization scheme these two noncommutative surfaces enter 
in a parallel fashion.  

For the quantum group $SO_q (4)$ the second and the third 
quadruplets  in the set (\ref{eq:quadruplet}) correspond to the
constraints originating from (\ref{eq:PMxx}), namely
\begin{eqnarray}
&&x_1 x_4 + q\, x_2 x_3 - q^{- 1}\, x_3 x_2 - x_4 x_1 = 0,
\nonumber\\   
&&x_1 x_4 - q^{- 1}\, x_2 x_3 + q\, x_3 x_2 - x_4 x_1 = 0.
\label{eq:x1x4}
\end{eqnarray} 
The consistency with the commutative limit is, therefore, 
maintained. The other two quadruplets in the list 
(\ref{eq:quadruplet}) correspond to the action of the 
diagonalized projectors $P_{(+)}$ and $P_{(0)}$ on the tensor 
product space ${\hat M} ( x \otimes x)$, and thereby lead to 
the $q$-deformation of the surfaces
\begin{equation}
x_1 x_4 + x_2 x_3 = {\sf k}_1,\qquad
x_1 x_4 - x_2 x_3 = {\sf k}_2.
\label{eq:surfSO4}
\end{equation}
Changing the basis from $( x_1, x_4 )$ and $( x_2, x_3 )$ to 
$( x \pm i t )$ and $( y \pm i z )$ respectively, we obtain a 
$3$-sphere in the first case ( with ${\sf k}_1 \geq 0$ ) and a
non-compact surface in the second case, as obtained before in 
the second equation in (\ref{eq:NCsurface}). Suitably 
implementing the $q$-dependent star operation we may obtain the 
corresponding noncommutative deformations related to the first and 
the last quadruplets in (\ref{eq:quadruplet}). 
  
\bigskip

\noindent{\Large{\bf Appendix A}}
\renewcommand{\theequation}{A.\arabic{equation}}
\setcounter{equation}{0}

\medskip 

\noindent Here we briefly indicate the derivations of the 
relations (\ref{eq:YYY}) and (\ref{eq:YXY}). We demonstrate this 
in the case of the quantum group $SO_q ( N )$, where the 
parameter $\varepsilon = 1$. To this end we define the operator 
\begin{equation}
K \equiv ( 1 + [ N - 1 ] ) P_{(0)} = \sum_{i, j = 1}^{N} 
q^{( \rho_i - \rho_j )} E_{i^{\prime}, j} \otimes 
E_{i, j^{\prime}},
\label{eq:KPEE}
\end{equation} 
where $i^{\prime} = N + 1 - i$. The relevant notations are 
explained in Sec. 4. Using the standard tensor structures
\begin{equation}
K_{12} = \sum_{i, j = 1}^{N}\, q^{( \rho_i - \rho_j )}\,
E_{i^{\prime}, j} \otimes E_{i, j^{\prime}} 
\otimes {\sf I}_{N \times N} \quad K_{23} = \sum_{i, j = 1}^{N}\, 
q^{( \rho_i - \rho_j )}\,{\sf I}_{N \times N} \otimes
E_{i^{\prime}, j} \otimes E_{i, j^{\prime}},
\label{eq:K12K13}
\end{equation}  
and the identity $E_{i, j^{\prime}} E_{k^{\prime}, l} = 
\delta_{j, k} E_{i, l}$, we obtain
\begin{equation}
K_{12} K_{23} = \sum_{i, j, k = 1}^{N} \,q^{( \rho_i - \rho_k )}
\,E_{i^{\prime}, j} \otimes  E_{i, k} \otimes E_{j, k^{\prime}}
\label{eq:Kprod}
\end{equation}  
and the following triple product rules
\begin{eqnarray}
K_{12} K_{23} K_{12} &=& \left( \sum_{i, j = 1}^{N} \,
q^{( \rho_i - \rho_j )}\,E_{i^{\prime}, j} \otimes  
E_{i, j^{\prime}} \right) \otimes \sum_{k = 1}^{N} E_{k, k}
\nonumber\\
\phantom{K_{12} K_{23} K_{12}} &=& K \otimes {\sf I}_{N \times N}
\equiv K_{12},\nonumber\\
K_{23} K_{12} K_{23} &=& \sum_{i = 1}^{N} E_{i, i} \otimes
\left( \sum_{j, k = 1}^{N}\, q^{( \rho_j - \rho_k )}\,
E_{j^{\prime}, k} \otimes E_{j, k^{\prime}} \right)\nonumber\\
\phantom{K_{23} K_{12} K_{23}} &=& {\sf I}_{N \times N} 
\otimes K \equiv K_{23}. 
\label{eq:Kident}
\end{eqnarray}
The triple product rules (\ref{eq:Kident}) and the defining 
property (\ref{eq:KPEE}) ensure that the following 
identities hold:
\begin{equation}
P_{(0)\;12} P_{(0)\;23} P_{(0)\;12} = ( 1 + [ N - 1 ] )^{-2} 
P_{(0)\;12},\quad
P_{(0)\;23} P_{(0)\;12} P_{(0)\;23} = ( 1 + [ N - 1 ] )^{-2} 
P_{(0)\;23}.
\label{eq:Pident}
\end{equation}
Using the definitions (\ref{eq:rhoval}) and (\ref{eq:PSp}) 
we may proceed analogously for the quantum group $Sp_q (N)$, 
where $N= 2 n$. The two cases can finally be unified after 
adopting the following definition
\begin{equation}
K \equiv \left( 1 + \varepsilon [ N - \varepsilon ] 
\right)\,P_{(0)}
\label{eq:KPboth}
\end{equation}
and then proceedinging as before. The final results may be 
summarized as 
\begin{equation}
P_{(0)\;12} P_{(0)\;23} P_{(0)\;12} = ( 1 + \varepsilon 
[ N - \varepsilon ] )^{-2} P_{(0)\;12},\quad
P_{(0)\;23} P_{(0)\;12} P_{(0)\;23} = ( 1 + \varepsilon
[ N - \varepsilon ] )^{-2} 
P_{(0)\;23}.
\label{eq:PSOSpid}
\end{equation}
Employing the definitions of $( Y_1, Y_2 )$ given in 
(\ref{eq:SOdef}), we now obtain (\ref{eq:YYY}). 

In order to prove the identity (\ref{eq:YXY}) we proceed as 
follows. Using the standard expressions for the braid 
generators ${\hat R}^{\pm 1}$ for the quantum groups 
$SO_q ( N )$ and $Sp_q ( N )$, and proceeding exactly 
analogously as before we obtain
\begin{equation}
P_{(0)\;12} {\hat R}^{\pm 1}_{23} P_{(0)\;12} = 
\frac{\varepsilon q^{\pm ( N - 1 - \varepsilon )}}
{ 1 + \varepsilon [ N - \varepsilon ]} P_{( 0 )\; 12},\quad 
P_{(0)\;23} {\hat R}^{\pm 1}_{12} P_{(0)\;23} = 
\frac{\varepsilon q^{\pm ( N - 1 - \varepsilon )}}
{ 1 + \varepsilon [ N - \varepsilon ]} P_{( 0 )\; 23}.
\label{eq:PRP}
\end{equation}
Expressing the braid matrices as in (\ref{eq:RRinvbrd}) we now 
use the relations (\ref{eq:PSOSpid}) and (\ref{eq:PRP}) 
to obtain
\begin{equation}
P_{(0)\;12} P_{(-)\;23} P_{(0)\;12} = \frac{\varepsilon
[ N - \varepsilon ] \left( [2] + 
\varepsilon [ N - 1 - \varepsilon ] \right)}
{ [2] \left( 1 + \varepsilon [ N - \varepsilon ] \right)^2}
P_{(0)\;12}.
\label{eq:P0PMP0}
\end{equation}  
The above result also holds after an exchange of the tensor 
indices: $(12) \rightleftharpoons (23)$. The equation 
(\ref{eq:P0PMP0}), in conjunction with the definitions 
(\ref{eq:Adef}) and (\ref{eq:SOdef}), now produces the 
identity (\ref{eq:YXY}).  

\bigskip

\noindent{\Large{\bf Appendix B}}
\renewcommand{\theequation}{B.\arabic{equation}}
\setcounter{equation}{0}

\medskip

\noindent  The correspondence between our result for 
Baxterization and that of \cite{I95} being quite simple in the 
context of the quantum group $GL_q (N)$, we discuss below the 
results for the quantum groups $SO_q (N)$ and $Sp_q (N)$. In 
Sec. 3.9 of \cite{I95} Isaev starts Baxterization of the braid 
matrices with the parametrization 
\begin{equation}
{\hat{\sf R}} ( x ) = c(x)\,\left( {\sf I} + a(x)\,
{\hat{\sf R}} + b(x)\,K \right),
\label{eq:RIsa}
\end{equation}
where the matrices ${\hat{\sf R}}$ and $K$ are given in 
(\ref{eq:BCDspcde}) and (\ref{eq:KPboth}) respectively. 
Substituting these results we may rewrite
\begin{equation}
{\hat{\sf R}} = c(x)\,( 1 + q a(x) )\, {\hat R} (x),
\label{eq:RIexp}
\end{equation}
where
\begin{equation}
{\hat R} (x) = {\sf I} - \frac{[2] a(x)}{1 + q a(x)}\,P_{(-)} +
\frac{\left( 1 + \varepsilon [ N - \varepsilon ] \right)\,b(x) - 
q \left( 1 - \varepsilon \,q^{- ( N + 1 - \varepsilon )} \right)
\,a(x)}{1 + q a(x)}\,P_{(0)}.
\label{eq:Icomp}
\end{equation}
Comparing this with our starting point (\ref{eq:BCDRvw}) we 
obtain the relations
\begin{equation}         
a(x) = - \frac{v(x)}{[2] + q v(x)},\,\,\,
b(x)= \frac{[2] ( w(x) - f_{+}\,v(x) )}{( [2] + q v(x) )\,
( 1 + \varepsilon [ N - \varepsilon ] )},\,\,\,
c(x) = [2]^{-1} ( [2] + q v(x) ),
\label{eq:Ifunction}
\end{equation}  
where the parameter $f_{+}$ has been defined in 
(\ref{eq:wvrat}). In the present work we have preferred the 
parametrization in (\ref{eq:BCDRvw}) as it assigns the key 
roles to the two projectors $P_{(-)}$ and $P_{(0)}$, leading 
to the systematic reduction of the trilinear forms presented in 
(\ref{eq:tlreduce}). As emphasized earlier this permitted us to 
display MBE and Baxterization as complementary facets of the 
same generalized ${\hat R} ( v, w )$ matrix. As for solutions of 
the braid equation obtained in the present work we note that the 
special case of the spectral variables $v = 0, w \neq 0$ has not 
been discussed in \cite{I95}. We have devoted Sec. 4 to study this 
remarkable new class of solutions. 

Our formalism led through (\ref{eq:xprime}) to the significant 
parallel structures of the Baxterized functions $v(x)$ and $w(x)$
given in (\ref{eq:vwx}). This enabled us to obtain the solutions 
{\em completely} by solving the two simplest equations, before
verifying that the more complex functional equation 
( corresponding to ${\sf a}_2$ in (\ref{eq:aabb}) ) is indeed 
consistent with them. In \cite{I95} the relatively complicated 
equation (Eq. (3.9.7) of \cite{I95}) had to be used to fix the 
two possible values of a parameter (Eq. (3.9.15) of \cite{I95}). 
The attractive forms of additive Baxterization, as evident in 
(\ref{eq:vthta})  and (\ref{eq:wSO1})-(\ref{eq:wSp2}) for the 
functions $v ( \theta )$ and $w ( \theta )$ respectively, are 
also direct consequences of our formalism.     
   
\noindent {\bf Acknowledgements}
\bigskip

\noindent 
{}
One of us (AC) thanks C. de Calan, J. Lascoux and J. Madore for 
discussions. The other author (RC) is partially supported by the
grant DAE/2001/37/12/BRNS, Government of India.

\bibliographystyle{amsplain}

\end{document}